\documentclass[12pt]{amsart}
\numberwithin{equation}{section}

\def\dim{\text{dim}\,}

\def\max{\text{max}\,}

\newcommand\trace{\operatorname{trace}}

\newcommand\x{\mathbf x}
\newcommand\y{\mathbf y}

\newcommand\C{\mathbf C}
\newcommand\F{{\mathbf F}}

\def\R{\mathbf R}

\newcommand\Z{\mathbf Z}

\def\1{\{1\}}
\def\2{\{2\}}
\def\3{\{3\}}
\def\4{\{4\}}

\newtheorem{theorem}{Theorem}
\newtheorem{lemma}[theorem]{Lemma}

\newtheorem{corollary}[theorem]{Corollary}
\newtheorem{proposition}[theorem]{Proposition}

\theoremstyle{definition}
\newtheorem{definition}[theorem]{Definition}
\newtheorem{example}[theorem]{Example}
\newtheorem{question}[theorem]{Question}
\theoremstyle{remark}
\newtheorem{remark}[theorem]{Remark}

\begin{document}

\title{On Kuiper's Conjecture}

\author[Cecil]{Thomas E. Cecil}
\thanks{The first author was partially supported by NSF Grant
No. DMS-0405529} 
\address{Department of Mathematics and Computer Science \\ College 
of the Holy Cross \\ Worcester, Massachusetts 01610-2395} 
\email{cecil@mathcs.holycross.edu}
\author[Chi]{Quo-Shin Chi}
\thanks{The second author was partially supported by NSF Grant
No. DMS-0103838} 
\address{Department of Mathematics \\ Campus Box 1146 \\ Washington
University \\ St. Louis, Missouri 63130} 
\email{chi@math.wustl.edu}
\author[Jensen]{Gary R. Jensen}
\thanks{The third author was partially supported by NSF Grant
No. DMS-0604236} 
\address{Department of Mathematics \\ Campus Box 1146 \\ Washington
University \\ St. Louis, Missouri 63130} 
\email{gary@math.wustl.edu}

\date{\today}


\begin{abstract} We prove that any connected proper Dupin hypersurface
in $\R^n$ is analytic algebraic and is an open subset of a connected
component of an irreducible algebraic set.
We prove the same
result for any connected non-proper Dupin hypersurface in $\R^n$ that satisfies a certain
finiteness condition.  
Hence any taut submanifold $M$ in $\R^n$, whose tube $M_\epsilon$ satisfies this finiteness condition,
is analytic algebraic and is a connected component of an irreducible algebraic set. 
In particular, we prove that every taut submanifold
of dimension $m \leq 4$ is algebraic.
\end{abstract}

\maketitle

\section{Introduction}\label{intro}
An embedding $f$ of a compact, connected manifold $M$ into Euclidean space $\R^n$ is
\textit{taut} if there is a field $\F$ such that every nondegenerate (Morse) Euclidean distance function, 
\[
L_p:M \to \R, \quad L_p(x) = d(f(x),p)^2, \quad p \in \R^n,
\]
has $\beta (M,\F)$ critical points on $M$, where $\beta (M,\F)$ is the sum of the $\F$-Betti numbers
of $M$.  That is, $L_p$ is a perfect Morse function on $M$. 

We can also consider taut embeddings into $S^n$, in which
case we use spherical distance functions instead of Euclidean distance functions.
Tautness is preserved by stereographic projection,
and so the theories in the two ambient spaces are essentially the same, and we will use whichever is most
appropriate for the particular question at hand.

Examples of taut submanifolds
are metric spheres, cyclides of Dupin (Banchoff \cite{Ban}), isoparametric hypersurfaces in spheres 
(Cecil-Ryan \cite{CecilRyan3})
isoparametric submanifolds of higher codimension (Terng \cite{Terng}),
and compact proper Dupin hypersurfaces (Thorbergsson \cite{Th1}).

In a paper published in 1984, Kuiper \cite{Kuiper2} asked whether all taut submanifolds of $\R^n$ must be real algebraic.
The affirmative answer to this question is now known as \textit{Kuiper's Conjecture}.  In the 1980's,
Kuiper's Conjecture was widely thought to be true, but a proof was never published. In particular, 
in 1984 Ulrich Pinkall sent
a sketch~\cite{P1} of a proof that proper Dupin hypersurfaces are algebraic in a letter to Tom Cecil, and we wish to
thank Professor Pinkall for allowing us to use his approach to the problem in our work in this paper. 

Pinkall's sketch led us to the local parametrization~(\ref{A18}) of a general proper Dupin hypersurface.  
We then used ideas from real algebraic geometry to show
that a connected proper Dupin hypersurface is contained in a connected component of an 
irreducible algebraic subset of $\R^n$.
There are still issues to be resolved to prove Kuiper's Conjecture, however, because arbitrary taut submanifolds are 
Dupin, but not necessarily proper Dupin.
 
In this paper, we prove Kuiper's Conjecture for manifolds of dimension $m \leq 4$, and we provide a criterion 
for algebraicity that may be useful in proving the conjecture in its entirety.  We now discuss the problem
in more detail.\\

Let $M$ be an immersed hypersurface in ${\R}^n$
or the unit sphere $S^n$ in ${\R}^{n+1}$. A {\em curvature surface} 
of $M$ is a smooth connected submanifold $S$
such that for each point $x \in S$, the tangent space
$T_xS$ is equal to a principal space of the shape operator
$A$ of $M$ at $x$.  The hypersurface $M$
is said to be \textit{Dupin} if it satisfies the condition

\begin{enumerate}
\item[(a)] along each curvature surface, the corresponding principal
curvature is constant.
\end{enumerate}
The hypersurface $M$ is called {\em proper Dupin} if, in addition 
to condition (a), it also satisfies the condition

\begin{enumerate}
\item[(b)] the number $g$ of distinct principal curvatures is constant
on $M$.
\end{enumerate}
Since both of these conditions are invariant under conformal
transformations, stereographic projection allows us to consider our
hypersurfaces in $\R^n$ or $S^n$, whichever is more convenient.

A primary result in this paper, contained in Theorem~\ref{AT2} and
Corollary~\ref{Main}, is a proof that any connected proper Dupin
hypersurface $M$ in $\R^n$ is analytic algebraic and thus is an open
subset of a connected component of an irreducible algebraic set in
$\R^n$.  If $M$ is complete, then it equals the connected component.

Pinkall \cite{Pinkall2} extended the notions of Dupin 
and proper Dupin to submanifolds $M$ of $\R^n$ of codimension greater than one.  In that
case, the principal curvatures are defined on the unit normal bundle $UN(M)$
of $M$ in $\R^n$.  Using Pinkall's definition, one can show that a submanifold
$M$ of $\R^n$ of codimension greater than one is proper Dupin if and only if
a tube $W^{n-1}$ over $M$ of sufficiently small radius is a proper Dupin hypersurface.  Since
algebraicity is also preserved by the tube construction, Theorem~\ref{AT2}
also implies that a proper Dupin submanifold of arbitrary codimension in $\R^n$
is algebraic.

These results can also be extended to non-proper Dupin hypersurfaces that satisfy what we call the local finiteness property (see Definition~\ref{finite}), which, for example, is implied by the condition
that any point not in the open dense subset ${\mathcal G}$ on which the multiplicities of the principal curvatures are locally constant has a
neighborhood that intersects ${\mathcal G}$ in only finitely many components. To our knowledge, 
this local finiteness property 
is satisfied by all known examples of Dupin hypersurfaces.

Thorbergsson \cite{Th1} proved that a compact proper Dupin hypersurface embedded in $\R^n$ is taut. 
Thus, our Theorem~\ref{AT2} shows that these taut hypersurfaces are algebraic.
Conversely, every taut submanifold   
in $\R^n$ is Dupin (see Miyaoka~\cite{Mi4} for hypersurfaces
and Pinkall~\cite{Pinkall2} in general), but not necessarily proper
Dupin (see Pinkall \cite{Pinkall}). In \cite{Pinkall2}, Pinkall also showed that
the tube $M_{\epsilon}$ of radius $\epsilon$ around an
embedded compact submanifold $M$ of $\R^n$ is $\Z_2$-taut if and only if $M$
is $\Z_2$-taut.
(The field $\F = \Z_2$ has proven to work well in the theory of taut submanifolds, and we will
use it exclusively in this paper.)  

In Theorem~\ref{Kuiper}, we prove: If $M$ is a compact taut submanifold of
$\R^n$ such that the tube $M_\epsilon$ over $M$, which is Dupin,
satisfies the local finiteness property, 
then $M$ is an analytic submanifold and a connected
component of an irreducible
algebraic subset of $\R^n$.
The proof of this theorem shows that the conclusion holds for any analytic
taut submanifold $M$, without the need to assume the local finiteness property on $M_\epsilon$.

In Section \ref{sec:6}, we develop the theory of Alexander cohomology, in
conjunction with the notion of ends, for noncompact manifolds needed in
Section \ref{sec:7} 
to prove that any taut submanifold of dimension 3 or 4 satisfies the local finiteness property given in
Theorem~\ref{Kuiper} and is therefore algebraic.  A key step in this approach is to show that the complement ${\mathcal Z}$
of the focal set ${\mathcal F}$ of a taut hypersurface $M$ in $S^n$ is connected, and that the normal exponential map $E$
from the normal bundle of $M$ to $S^n$ restricts to a finite covering map on $E^{-1}({\mathcal Z})$.
These results, together with the classifications of taut submanifolds of dimensions one and two due to Banchoff \cite{Ban},
show that Kuiper's Conjecture is true for taut submanifolds
of dimension $m \leq 4$ (Theorem~\ref{dim<5}). This approach may have applications
to higher dimensional taut submanifolds, although we have not been able to extend it further at this time.

The paper is organized as follows.
Section \ref{sec:1} outlines the facts needed about algebraic and
semi-algebraic subsets of $\R^n$ and explains that Nash functions are
the same as analytic algebraic functions.  
Our principal algebraic
tool is contained in Lemma~\ref{lemma0}, which states that if a
connected analytic submanifold $M \subset \R^n$ contains a connected
open subset that is a semi-algebraic subset of $\R^n$, then $M$ is an
open subset of a connected component of an irreducible algebraic set.

Section \ref{sec:2} outlines how to calculate the center and radius of
a spherical curvature surface in terms of the mean curvature normal.

In  Section \ref{sec:3}, we present the facts needed about jet spaces and show
that the curvature surface through a point for a principal curvature
of constant multiplicity is determined by the 3-jet
at the point of the embedding.

The proof of Theorem~\ref{AT2} showing that proper Dupin hypersurfaces are algebraic is contained in Section \ref{sec:4}.

In Section \ref{sec:5}, we define the local finiteness property and
extend the results of Section \ref{sec:4} to any connected Dupin hypersurface
possessing this finiteness property.  We then prove Theorem~\ref{Kuiper} concerning taut submanifolds mentioned above.

Finally, in Section \ref{sec:6}, we develop the theory of Alexander cohomology and ends for noncompact manifolds, and in 
Section \ref{sec:7}, we prove that any taut submanifold of dimension $m \leq 4$ in $\R^n$ is algebraic.

\section{Algebraic Preliminaries}\label{sec:1}
In this section we briefly 
review the material we need in real algebraic geometry,
referring for details to the book by Bochnak, et al.~\cite{RAG}.

Let $\R[X_1,\dots,X_n]$ denote the ring of all real polynomials in $n$
variables.  If $B$ is any subset of this ring, we denote the
\textit{set of zeros} of $B$ by
\[
\mathcal Z(B) = \{x \in \R^n : \forall f \in B, \; f(x)=0\}.
\]
\begin{definition}\label{algebraic}
An \textit{algebraic subset} of $\R^n$ is the set of zeros of some $B
\subset \R[X_1,\dots,X_n]$.  
\end{definition}
Given any subset $S$ of $\R^n$, the set
of all polynomials vanishing on $S$, denoted
\[
\mathcal I (S) = \{f \in \R[X_1,\dots,X_n]: \forall x \in S, \;
f(x)=0\},
\]
is an ideal of $\R[X_1,\dots,X_n]$.  Ideals of the ring of real
polynomials are finitely generated, which implies that for any
algebraic set $V$ of $\R^n$, there exists a polynomial $f \in
\R[X_1,\dots,X_n]$ such that $V = \mathcal Z(f)$.  In fact, if
$f_1,\dots,f_m$ generate $\mathcal I(V)$, take
\[
f = f_1^2 +\dots + f_m^2.
\]

\begin{definition}\label{Def2.1}
A \textit{semi-algebraic subset} of $\R^{n}$ is one which is a finite
union 
of sets of the form
$$
\cap_{j} \{x\in \R^{n}: F_j(x) * 0\},
$$
where * is either $<$ or $=$, $F_j\in \R[X_1,\dots,X_n]$, and the
intersection is finite. ($F_j<0$ is just 
$-F_j>0$. 
So there is no need to introduce $>$ in the definition.)
\end{definition}

All algebraic sets are certainly semi-algebraic, and $\R^{n}$
itself is given by 
the null relation. Another simple example popular in calculus is  
$$
\{(x,y): y-1+x^2\leq 0\}\cap\{(x,y):x^2-y\leq 0\},
$$
the area enclosed by the two parabolas.

It follows from the definition that the complement of a semi-algebraic
set is 
semi-algebraic, and hence a semi-algebraic set take away another
semi-algebraic set leaves a semi-algebraic set.

The following result is proved in \cite[Theorem 2.21, p.\ 26]{RAG}.
\begin{proposition}\label{th0}
The projection $\pi:\R^n\rightarrow \R^k$ sending
$x\in\R^{n}$ to its first $k$ coordinates maps a
semi-algebraic set to a 
semi-algebraic set.
\end{proposition}

In particular, even though a linear projection does not map an algebraic set
to an algebraic set, in general, it does map it
to a semi-algebraic 
set. 

\begin{corollary}\label{Cor2.1}
The (topological) closure and interior of a semi-algebraic set are
semi-algebraic. 
\end{corollary}
A proof is given in \cite[Proposition 2.2.2, p.\ 27]{RAG}.

\begin{definition}\label{Def2.2}
A map $f:S\subset \R^{n}\rightarrow\R^{k}$
is \textit{semi-algebraic} if its graph in $\R^{n}\times\R^{k}$
is a semi-algebraic set.
\end{definition}
\begin{corollary}\label{Cor2.2}
The image of a semi-algebraic map
$f:S\subset \R^{n}\rightarrow\R^{k}$
is semi-algebraic, via the composition
$$
{\rm
graph}(f)\subset\R^{n}\times\R^{k}\rightarrow\R^{k},
$$
where the last map is the projection onto the second summand.
\end{corollary}


\begin{definition}\label{Def2.3}
A \textit{Nash function} is a $C^{\infty}$ semi-algebraic map from an open
semi-algebraic subset of $\R^{n}$ to $\R$.
\end{definition}

\begin{definition}\label{Def2.4} A real analytic function $f$ defined
on an open 
semi-algebraic subset $U$ of $\R^n$ is \textit{analytic algebraic} if it
is a solution of a polynomial equation on $U$ of the form,
\begin{equation}\label{Eq2.9}
a_{0}(x) f^s(x) + a_{1}(x) f^{s-1}(x) + \cdots + a_{s}(x) = 0,
\end{equation}
where $a_0(x)\neq 0,a_{1}(x), \cdots, a_s(x)$ are polynomials over
$\R^{n}$.
\end{definition}
A significant result in the subject is that these two concepts are
equivalent (see \cite[Proposition 8.1.8, p.\ 165]{RAG}).  
\begin{proposition}\label{char}
A function is Nash if and
only if it is analytic algebraic.
\end{proposition}
The following example will be useful to us later in the paper.
\begin{example}\label{Ex2.1} For any number $\epsilon$ satisfying
$0 < \epsilon < 1$, the open ball
\begin{equation}\label{Eq2.10}
B^n(\epsilon) = \{s = (s^1,\dots,s^n)\in \R^n : |s| < \epsilon \}
\end{equation}
is an open semi-algebraic subset of $\R^n$.  The function
\begin{equation}\label{Eq2.11}
s^0 = \sqrt{1-|s|^2}
\end{equation}
on $B^n(\epsilon)$ is analytic algebraic, since $(s^0(x))^2 + a_0(x) =
0$ on $B^n(\epsilon)$, where $a_0(x)$ is the polynomial $|s|^2-1$ on
$\R^n$. Partial derivatives of all orders of $s^0$ are analytic
algebraic.  In fact, an elementary calculation and induction argument
shows that if $D_i$ denotes the partial derivative with respect to
$s^i$, then
\[
D_{i_1\dots i_k} s^0 = \frac{a_k(s)}{(s^0)^m}
\]
where $a_k(s)$ is a polynomial on $\R^n$ and $m$ is a positive
integer.  Therefore, 
\[
(s^0)^{2m} (D_{i_1\dots i_k} s^0)^2-a_k(s)^2 = 0
\]
is an equation of the form~(\ref{Eq2.9}), since $(s^0)^2$ is a
polynomial on $\R^n$.
\end{example}

A slight generalization of the single-variable case in \cite[Proposition
2.9.1, p.\ 54]{RAG}, shows that the partial derivatives of any Nash
function are again Nash functions.  
\begin{proposition}\label{good}
Let $f$ be a Nash function defined on a connected semi-algebraic
open set $U$ in $\R^{n}$ .
Then a partial derivative of $f$ of any order with respect to 
$x=(x^{1},\cdots,x^{n})\in\R^{n}$ is again a Nash function on $U$.
\end{proposition}
\begin{definition}
Let $S$ be a semi-algebraic subset of $\R^{n}$. The
\textit{dimension} of $S$, denoted $\dim S$, is the dimension of the ring
$R = \R[x^1,\cdots,x^n]/\mathcal I(S)$,  which is 
the maximal length of chains of prime ideals of $R$.
\end{definition}
In \cite[Proposition 2.8.14, p.\ 54]{RAG} it is proved that
\begin{proposition}\label{dimension}
If $S$ is a semi-algebraic subset of $\R^n$ that is a $C^\infty$
submanifold of $\R^n$ of dimension $d$, then $\dim S = d$.
\end{proposition}

We denote the topological closure of $S$ by $\overline{S}^{\text{top}}$. 
Denote the Zariski closure of $S$ by $\overline{S}^{\text{zar}}$.  It is the
smallest algebraic set containing $S$, so
\[
\overline{S}^{\text{zar}} = \mathcal Z(\mathcal I(S)).
\]
The following plays a central role in our arguments to follow
(see \cite[Proposition 2.8.2, p.\ 50]{RAG}). 
\begin{proposition}\label{prop}
If $S$ is a semi-algebraic subset of $\R^{n}$, then
\[
\dim S=\dim(\overline{S}^{\text{top}})=\dim(\overline{S}^{\text{zar}}).
\]
\end{proposition}
\begin{definition}\label{irreducible}
An algebraic set $V$ in $\R^n$ is \textit{irreducible} if, whenever $V
= F_1 \cup F_2$, where $F_1$ and $F_2$ are algebraic sets, then
$V=F_1$ or $V = F_2$.
\end{definition}
The following facts are in \cite[Theorem 2.8.3, p.\ 50]{RAG}:
Every algebraic set $V$ is the union, in a unique way, of a finite
number of irreducible algebraic sets $V_1,\dots,V_p$, such that
$V_i$ is not a subset of $\cup_{j\neq i} V_j$
for $i = 1,\dots,p$.  The $V_i$ are the \textit{irreducible
components} of $V$ and
\[
\dim V = \max(\dim V_1,\dots,\dim V_p).
\]
An algebraic set $V \subset \R^n$ is irreducible if and only if the
ideal $\mathcal I(V)$ is prime.

\begin{remark}
Let $X$ be an irreducible real algebraic set in
$\R^{n}$. By~\cite{Hironaka}, $X$ can be successively blown up to 
reach its desingularization $X^{*}$. Let $\pi:X^{*}\to X$ be the
projection. Then the proper transform 
of $X$ in $X^{*}$, that is, the Zariski closure of the preimage of $X$
with its singular set removed, is a smooth  
algebraic set with disjoint smooth connected (algebraic) components
$X^{*}_{1},\ldots,X^{*}_{s}$, whose projections in $X$ are called the
{\em connected components} of $X$. We will show 
these connected components are semi-algebraic as follows.  

$X^{*}_{1},\ldots,X^{*}_{s}$ are exactly the irreducible components of
the proper transform of 
$X$ under $\pi$. The morphism $\pi$ is algebraic. Hence it maps the
algebraic sets  
$X^{*}_{i}$ to the connected components of $X$, which are therefore
semi-algebraic. Note that the proper transform  
itself is not affine, but we can always cover it by a finite number of
affine charts and then proceed with our arguments 
in each chart and take a finite union. 

That the connected components of $X$ are closed follows from the
fact that $\pi$ is a proper map.  

A real irreducible algebraic set may have several connected
components of varying dimensions. 
A good example is a degenerate torus given by
$$
y^2=x^2(x-1).
$$
The solution set is irreducible since the polynomial is. This
irreducible algebraic set 
consists of two connected components. One is the curve that opens
toward positive 
infinity, with vertex at $(1,0)$. The other is $\{(0,0)\}$. 

The name connected component is standard in real algebraic geometry,
but we should realize that these are not always the connected
components in the topological sense.
For example, the irreducible  
algebraic set known as 
the Cartan umbrella (see \cite[pp.\ 60-61]{RAG}).
$$
z(x^2+y^2)-x^3=0
$$
is topologically connected, but has two connected components:
the surface $z=x^3/(x^2+y^2)$ (including the origin) and the
$z$-axis. 

If an irreducible algebraic set $X$ in $\R^n$ contains a connected
analytic hypersurface $M$ in $\R^n$, then $M$ must be an open subset
of one of the connected components of $X$.  If $M$ is closed, then it
must be the whole connected component.
\end{remark}

We can now prove our principal algebraic tool.

\begin{lemma}\label{lemma0}
Let $M \subset\R^{n}$ be a connected analytic submanifold and
let $U\subset M$ be a connected open subset of $M$.  If $U$
is a semi-algebraic subset of $\R^{n}$, then $M$
is an open subset of a connected component of
the irreducible algebraic
set $\overline{U}^{\text{zar}}$.
\end{lemma}

\begin{proof} The Zariski closure $\overline{U}^{\text{zar}}$ of $U$ is
irreducible, because $\mathcal I(\overline{U}^{\text{zar}}) = \mathcal I(U)$
is a prime ideal.  In fact, suppose $f$ and $g$ are polynomials such
that $fg \in \mathcal I(U)$.  If $f \notin \mathcal I(U)$, then there
exists a point $p\in U$ such that $f(p) \neq 0$.  By continuity, there
exists an open subset $O$ of $U$, containing $p$, on which $f$ is
never zero.  But $fg$ is identically zero on $O$, so $g$ must be
identically zero on $O$.  Since $O$ is an open subset of $U$, and
since $g$ is an analytic function on $U$, it follows that $g$ is
identically zero on the connected set $U$, that is, $g \in \mathcal
I(U)$.  

By the same sort of argument, if $f \in \mathcal I(U)$, then $f$ is an
analytic function on $M$, identically zero on the open subset $U$, so
must be identically zero on $M$.  Hence $M \subset \mathcal Z(\mathcal
I(U)) = \overline{U}^{\text{zar}}$.  
\end{proof}

\begin{definition}\label{NashSub}
A semi-algebraic subset $M$ of $\R^{m}$ is a \textit{Nash submanifold} of
$\R^{m}$ of dimension 
$n$ if for every point $p$ of $M$, there exists a Nash diffeomorphism
$\psi$ from an open semi-algebraic  
neighborhood $U$ of the origin in $\R^{m}$ into an open
semi-algebraic neighborhood $V$ 
of $p$ in $\R^{m}$ such that $\psi(0)=p$ and
$\psi((\R^{n}\times\{0\})\cap U)= M\cap V$. Here, by a \textit{Nash  
diffeomorphism} $\psi$ we mean the coordinate
functions of $\psi$ and $\psi^{-1}$ are Nash functions.
\end{definition}
\begin{definition} Let $M$ be a Nash submanifold of
$\R^{m}$. A mapping $f:M\rightarrow\R$ is a  
\textit{Nash mapping} if it is semi-algebraic, and for every $\psi$ in the
preceding definition, $f\circ\psi|_{\R^{n}\cap U}$ 
is a Nash function.
\end{definition}
Semi-algebraic subsets of $\R^n$
are, in a sense, 
piecewise algebraically analytic. This is made precise in the
following result (see \cite[Proposition 2.9.0, p.\ 57]{RAG}).   
\begin{proposition}\label{th1}
Let $S$ be a semi-algebraic subset of $\R^{n}$. Then $S$ is the
disjoint union of a finite number of Nash submanifolds $M_i$ of
$\R^{n}$, 
each Nash diffeomorphic to an open cube $(0,1)^{\dim M_{i}}$.
\end{proposition}

The following lemma is a kind of converse to Lemma~\ref{lemma0}, in
that it implies that if a $C^\infty$ submanifold $U^d$ of $\R^n$ is
contained in an irreducible algebraic subset of dimension $d$ in
$\R^n$, then $U$ is a real analytic submanifold.
    
\begin{lemma}\label{keylemma}Let $X$ be a closed semi-algebraic subset
of 
$\R^{n}$ of dimension 
$d$. If there is an open set $U$ of $X$ such that
$U$ is a $C^\infty$ submanifold of dimension $d$ in $\R^{n}$,
then around each point $p$ of $U$ there is an open neighborhood $B
\subset U$ that is a
semi-algebraic subset of $\R^n$.
\end{lemma}
\begin{proof} By Proposition~\ref{th1}, $X$ is the finite disjoint union of
Nash submanifolds, each Nash diffeomorphic to an open
cube of some dimension. Let $M_{1}, \ldots,
M_{k}$ be those Nash submanifolds of dimension $d$ that intersect $U$. There
are such Nash submanifolds, for otherwise, $U$ would be contained in
a finite union of Nash submanifolds 
of smaller dimensions. Let $Y$ be the topological closure of the union of
$M_{1},\ldots,M_{k}$. Then $Y$ is semi-algebraic and $Y\subset X$.

$U$ is contained in $Y$, for otherwise, there would be a point $z\in U$
that is not in $Y$, so that $z$ stays some distance away from
$Y$. Therefore, 
there would be a small open ball $D$ around $z$ in $U$ disjoint from $Y$,
which would imply that $D$ is contained in Nash manifolds of smaller
dimensions, 
because $Y$ contains all the Nash submanifolds of dimension $d$ not disjoint
from $U$. This is absurd.

The part $Y^{*}$ where $Y$ is not in
$U$ is closed. Since $p\in U$, there is a nonzero shortest distance $\delta$
from $p$ to $Y^{*}$. Choose a Euclidean ball $B$ of radius $\delta'<\delta$
centered at $p$ such that $B\cap U$ is a submanifold of $\R^{n}$
diffeomorphic to a cube of dimension $d$, which is possible since $U$ is a
submanifold of $\R^{n}$.

We claim that $B\cap U=B\cap Y$. It is clear that $B\cap U\subset B\cap Y$
since $U\subset Y$. Suppose $w\in B\cap Y$. Then $w$ must be in $U$, because otherwise
$w\in Y^{*}$, so that the distance from $p$ to $w$ would be $> \delta'$ and
so $w$ would not be in $B$, a contradiction. Hence, $w\in B\cap U$, proving the claim.

Now that $B\cap U=B\cap Y$, we see that $B\cap U$ on the one hand is
diffeomorphic to an open
cube of dimension $d$ contained in $U$, and on the other hand it is
semi-algebraic since so are $B$ and $Y$.
\end{proof}

\section{Mean curvature normal}\label{sec:2}

\begin{definition}\label{umbilic} Following
Cecil-Ryan~\cite[p.\ 140]{CecilRyan}, we say that a submanifold
$\x:V^m \to 
\R^n$, for any $m<n$, is \textit{umbilic} if there exists a linear map
$\omega_{(x)}: T_x^\perp V \to \R$, for every $x \in V$, such that for
every $Z \in T_x^\perp V$, the shape operator
\begin{equation}\label{B1}
B_Z :T_xV \to T_xV \quad \mbox{ is } \quad B_Z = \omega(Z)I.
\end{equation}
\end{definition}
Let $\x,e_i,e_\alpha$ be a \textit{first order} frame field on an open set
$U\subset V$ along $\x$, which means that the vector fields
$e_i$, for $i = 1,\dots,m$, are tangent to $\x(V)$ and the $e_\alpha$, for
$\alpha = m+1,\dots,n$, are normal to $\x(V)$ at each point of $U$.  Then 
\begin{equation}\label{B2}
d\x = \sum_i \theta^i e_i, \quad de_i = \sum_j \omega^j_i e_j +
\sum_\alpha \omega^\alpha_i e_\alpha, \quad de_\alpha = \sum_j \omega^j_\alpha e_j +
\sum_\beta \omega^\beta_\alpha e_\beta
\end{equation}
where the 1-forms $\omega^a_b = - \omega^b_a$, for $a,b = 1,\dots,n$,
and
\begin{equation}\label{B3}
\omega^\alpha_i = \sum_{i,j} h^\alpha_{ij}\theta^j
\end{equation}
for functions $h^\alpha_{ij} = h^\alpha_{ji}$.  Then 
\begin{equation}\label{B4}
B_{e_\alpha} = \sum_i\omega^\alpha_i e_i = \sum_{i,j} h^\alpha_{ij}\theta^j e_i.
\end{equation}
A normal vector is given by $Z = \sum_\alpha z^\alpha e_\alpha$, so
\[ 
B_Z = \sum_{\alpha}z^\alpha B_{e_\alpha} = \sum_{i,j,\alpha} z^\alpha
h^\alpha_{ij}\theta^j e_i.
\]  

\begin{definition}\label{mean} The \textit{mean curvature normal}
vector field of $\x$ is the normal vector field on $V$
\[  
H = \sum_\alpha (\frac 1m \trace B_{e_\alpha})e_\alpha = \frac 1m \sum_\alpha
\left( \sum_i h^\alpha_{ii}\right) e_\alpha,
\]  
which is independent of the choice of orthonormal frame $e_\alpha$ in the
normal bundle.
\end{definition}
We see that $\x$ is
umbilic if and only if $B_{e_\alpha} = h^\alpha I$, for some function
$h^\alpha$, for all $\alpha$; that is,
\[  
\omega^\alpha_i = h^\alpha \theta^i,
\]  
for all $i,\alpha$.  When $\x$ is umbilic, then the mean curvature normal
is
\[  
H = \sum_\alpha h^\alpha e_\alpha.
\]  

Let $\x:M^{n-1} \to \R^n$ be a connected immersed hypersurface.  Let
$\xi$ be a field of unit normals to $\x(M)$ defined on an open subset
$U$ of $M$, and let $A$ denote the shape operator corresponding to
$\xi$.  If the principal curvatures are ordered
\[
\kappa_1(p) \leq \dots \leq \kappa_{n-1}(p)
\]
for each $p\in U$, then the $\kappa_i$ are continuous functions on $U$
(see Ryan \cite[Lemma 2.1, p.271]{Ryan}).  If a continuous principal
curvature function $\kappa$ has constant multiplicity $m$ on $U$, then
$\kappa$ is smooth on $U$, as is its $m$-dimensional principal
distribution $T_\kappa$ of eigenvectors of $A$ corresponding to the
eigenvalue $\kappa$ (see Nomizu \cite{Nomizu} or Singley
\cite{Singley}).  Furthermore, $T_\kappa$ is integrable and we will
refer to $T_\kappa$ as the \textit{principal foliation} corresponding
to the principal curvature $\kappa$.

\begin{proposition}\label{BC1}  Suppose $\x: M^{n-1} \to \R^n$ is a
connected 
immersed hypersurface with a principal curvature function $\kappa$ of
constant 
multiplicity $m$.  If $m=1$, suppose further that $\kappa$ is constant
on any leaf $V^1$ of the principal foliation
$T_\kappa$ of $M$.  For any leaf $V^m$ of $T_\kappa$,

i). If $\kappa = 0$ and $d\kappa = 0$ at each point
of $V$, then $\x(V)$ is an open subset of an
$m$-plane in $\R^n$. 

ii). If $\kappa$ or $d\kappa$ is nonzero at some
point of $V$, then $\x(V)$ is an open subset of an 
$m$-sphere in $\R^n$ of
radius $1/|H|$ and center at 
\begin{equation}\label{B75aa}
\x(p_0) + H(p_0)/|H|^2,
\end{equation}
where $p_0$ is any point in $V$ and
$H$ is the mean curvature normal of $\x(V^m) \subset \R^n$.
\end{proposition}

\begin{proof}  See \cite[Theorem 4.5, pp.\ 140-141 and Theorem 4.8,
pp.\ 145-147]{CecilRyan}.
\end{proof}

We shall need an explicit formula for the mean curvature vector $H$ of
$V$ in Proposition~\ref{BC1}.
Let $\x, e_i, e_\alpha, e_n$ be a first order frame field
along $\x$, for which $e_n$ is a unit normal vector field, each vector
field $e_i \in T_\kappa$, 
where $i = 1,\dots,m$, and each vector field $e_\alpha \in T^\perp_\kappa
\subset 
T_\x M$, where $\alpha = m+1,\dots,n-1$.  Then $d\x = \sum_i \theta^i e_i
+ \sum_\alpha \theta^\alpha e_\alpha$ and $d e_n = \sum_{a=1}^{n-1}
\omega^a_n e_a$ and the shape operator $B_{e_n}: T_xM \to T_xM$
satisfies $B_{e_n} e_i = \kappa e_i$, for all $i$, and $B_{e_n} e_\alpha
\in \mbox{span}\, \{e_\beta\}$, for all $\alpha$, since $B_{e_n}$ is self-adjoint.
But
\begin{equation}\label{19}
B_{e_n} = \sum_{a=1}^{n-1} \omega^n_a e_a = \sum_i \omega^n_i e_i +
\sum_\alpha \omega^n_\alpha e_\alpha,
\end{equation}
where $\omega^n_a = \sum_{b=1}^{n-1} h_{ab} \theta^b$.  Therefore,
\begin{equation}\label{B20}
\omega^n_i = \kappa \theta^i, \quad \omega^n_\alpha = \sum_\beta h_{\alpha\beta}
\theta^\beta,
\end{equation}
for all $i,\beta$; that is, $h_{i\alpha} = h_{\alpha i} = 0$, for all $i,\alpha$.  In
addition, the $(n-m-1) \times (n-m-1)$ symmetric matrix
of functions
\begin{equation}\label{B21}
c = (c_{\alpha\beta}) =(h_{\alpha\beta} - \kappa \delta_{\alpha\beta}) 
\end{equation}
is non-singular at every point of $V$, so has inverse
\begin{equation}\label{B28b}
c^{-1} =
(c^{\alpha\beta}).
\end{equation}
Set
\begin{equation}\label{derivative}
d\kappa = \sum_i \kappa_i \theta^i + \sum_\alpha \kappa_\alpha \theta^\alpha.
\end{equation}
Under the hypothesis of Proposition~\ref{BC1}, the $\kappa_i=0$.  The
$\kappa_\alpha$ are the derivatives of $\kappa$ normal to the leaves of
$T_\kappa$.  At any point of a leaf $V^m$ of $T_\kappa$,
the mean curvature normal of $V$ at the point is
\begin{equation}\label{B75}
H = \kappa e_n - \sum_{\alpha,\beta} c^{\alpha\beta} \kappa_\beta\, e_\alpha.
\end{equation}

\section{Jets}\label{sec:3}

\begin{definition} Smooth maps $f,g: M^m \to N^n$ have the same
  $k$-jet at $p \in M$ means $f(p) = g(p) = q$, and for any smooth
  curve $\gamma: \R \to M$ such that $\gamma(0)=p$ and any smooth
  function $F:N \to \R$ such that $F(q)=0$, the derivatives
\begin{equation}\label{A1}
(F\circ f \circ \gamma)^{(j)} (0) = (F\circ g \circ \gamma)^{(j)}(0)
\end{equation}
for all $j \leq k$.  This defines an equivalence relation on the set
of smooth maps $C^\infty(M,N)$.  If $f \in C^\infty(M,N)$, denote its
$k$-jet at $p\in M$ by $j^k_p(f)$.
\end{definition}

\begin{proposition}\label{AT1}  If $f,g \in C^\infty(M,N)$, then
$j^k_p(f) = 
  j^k_p(g)$ if and only if for any local coordinate charts $U, x =
  (x^1,\dots,x^m)$ about $p$ in $M$ and $V, y=(y^1,\dots,y^n)$ about
  $f(p)=g(p)=q$ in $N$, all partial derivatives 
\begin{equation}\label{A2}
D_{i_1\dots i_j} f^\alpha(p) = D_{i_1\dots i_j} g^\alpha(p), \quad
\mbox{ for all $j \leq k$, $\alpha = 1,\dots,n$},
\end{equation}
where $f^\alpha = y^\alpha \circ f$ and $g^\alpha = y^\alpha \circ g$.
\end{proposition}

\begin{proof} See \cite[pp.\ 20-22]{BCGGG}.
\end{proof}

If $(p,q) \in M \times N$, let 
\begin{equation}\label{AA3}
J^k_{p,q}(M,N) = \{ j^k_p(f): f \in C^\infty(M,N), \; f(p)=q\},
\end{equation}
and let
\begin{equation}\label{A4}
J^k(M,N) = \cup_{(p,q) \in M\times N} J^k_{p,q}(M,N)
\end{equation}
for each whole number $k$.  This space has a natural $C^\infty$
structure with $C^\infty$ projections
\begin{equation}\label{A5}
\begin{array}{rlrl}
J^k(M,N) & \to M, &\quad J^k(M,N) &\to N \\
j^k_p(f) &\mapsto p &\quad j^k_p(f) &\mapsto f(p).
\end{array}
\end{equation}
A map $f \in C^\infty(M,N)$ defines a $C^\infty$ section
\begin{equation}\label{A6}
\begin{array}{rl}
j^k(f): M &\to J^k(M,N) \\
p &\mapsto j^k_p(f).
\end{array}
\end{equation}

\begin{lemma}\label{AL1} For any $f \in C^\infty(M,N)$, and $p\in M$,
the jet $j^{k+l}_p(f) \in J^{k+l}(M,N)$ can be identified with
$j^k_p(j^l(f)) \in J^k(M,J^l(M,N))$.
\end{lemma}

\begin{proof} Local coordinates in $M$ and $N$ determine local
coordinates in $J^l(M,N)$ (see \cite{BCGGG}), and these give the
local expression
\begin{equation}\label{A7}
j^l(f) = (p,f,D_{i_1\dots i_j}f), \quad j \leq l.
\end{equation}
Then $j^k_p(j^l(f))$ is given by the first $k$ partial derivatives of
these component functions at $p$, which are, ignoring redundancies,
the first $k+l$ partial derivatives of $f$ at $p$.
\end{proof}
It seems justified to abuse notation to express this lemma as
\begin{equation}\label{A7a}
j^{k+l}_p(f) = j^k_p(j^l(f)),
\end{equation}
for any $f \in C^\infty(M,N)$, and $p\in M$.

\begin{remark}\label{AR1}
Let $\y:V^m \to \R^n$ be an immersed submanifold.  Any first order
frame field 
\begin{equation}\label{A7aa}
\y, e_i, e_\alpha
\end{equation}
along $\y$, where $e_i$ are tangent, $i =
1,\dots,m$, and 
$e_\alpha$ are normal, $\alpha = m+1,\dots n$, is determined by the 1-jet 
$j^1_p(\y)$ at each point $p\in V$.  Now
\begin{equation}\label{A7b}
d\y = \sum_i \theta^i e_i, \quad de_a = \sum_{b=1}^n \omega^b_a e_b,
\; a = 1,\dots,n,
\end{equation}
so the 1-forms $\theta^i$ are determined by the
1-jet of $\y$ at each point and  the 1-forms $\omega^b_a$ are
determined by the 2-jet of $\y$ at each point, since each $e_a$ is
determined by the 1-jet of $\y$.  Then
\begin{equation}\label{A7c}
\omega^\alpha_i = \sum_j h^\alpha_{ij}\theta^j,
\end{equation}
where the functions $h^\alpha_{ij} = h^\alpha_{ji}$ are determined by the
2-jet of $\y$ at each point.  Hence, the shape operators~(\ref{B4})
are determined by the 2-jet of $\y$ at each point, and thus the shape
operator $B_Z$ is determined by the 2-jet of $\y$, for any normal
vector $Z = \sum_\alpha z^\alpha e_\alpha$.

In terms of our first order frame field~(\ref{A7aa}) on $U \subset V$,
the characteristic polynomial of the shape operator $B_{e_\alpha}$
of~(\ref{B4}) is
\begin{equation}\label{A7d}
F^\alpha(p,z) = \det (h^\alpha(p) - zI),
\end{equation}
for $p \in U$, $z \in \C$, where $h^\alpha(p) = (h^\alpha_{ij}(p))$ is the
$m\times 
m$ symmetric matrix defined in~(\ref{A7c}).  Thus, $F(p,z)$ is a
polynomial function of $h^\alpha_{ij}(p)$ and $z$, and is determined by
the 2-jet of $y$.  The eigenvalues of $B_{e_\alpha}$ are necessarily real,
so can be arranged as
\begin{equation}\label{A7e}
\kappa^\alpha_1(p) \leq \dots \leq \kappa^\alpha_m(p),
\end{equation}
for each $p \in U$.  These are the roots of $F^\alpha(p,z)$.  They are
continuous functions on $U$ (see Ryan \cite[Lemma 2.1, p.\ 271]{Ryan}).
A slight variation on Ryan's argument shows that if
$\kappa^\alpha_i$ has constant multiplicity $m_i$ on $U$, for some $i$,
then it 
is a rational function of the entries of $h^\alpha$ on $U$.  In
particular, $\kappa^\alpha_i$ is determined by the 2-jet of $\y$ at each
point, so its derivative $d\kappa^\alpha_i$ is determined by the 3-jet of
$y$.

Suppose now that $\x:M^{n-1} \to \R^n$ is a hypersurface, so that we
can drop the superscript $\alpha$ in the above notation.  Suppose that
$\kappa$ is a principal curvature of constant multiplicity $m$ on
$U \subset M$.  Then the
functions $c_{\alpha\beta}$ in~(\ref{B21}) and $c^{\alpha\beta}$ of~(\ref{B28b})
depend rationally on the 2-jet of
$\x$ at each point and the normal derivatives $d\kappa(e_\alpha) =
\kappa_\alpha$ depend rationally on the 3-jet of $\x$ at each point.
Therefore, the mean curvature vector $H$ in~(\ref{B75}) depends
rationally on
the 3-jet of $\x$ at each point.
\end{remark}

\section{Proper Dupin hypersurfaces are algebraic}\label{sec:4}

We are ready to prove that a proper Dupin hypersurface $f: M^{n-1} \to
\R^n$ is 
analytic algebraic in the sense that its differentiable structure has an
analytic subatlas with respect to which the immersion into $\R^n$ is
analytic algebraic.   
It follows that each point of $M$ is
contained in an open neighborhood of $M$ that is a semi-algebraic
subset of $\R^n$.  From the analyticity and the semi-algebraic
neighborhood  we can conclude that if $M$ is connected,
then it is an open subset of an
irreducible algebraic set in $\R^n$ of dimension $n-1$.

Suppose $f$ has  $g$
principal curvatures $\kappa_1,\dots,\kappa_g$,
with multiplicities $m_1,\dots,m_g$.  Let $\{1\},\dots,\{g\}$ be a
partition of the set $\{1,\dots,n-1\}$ into disjoint subsets for which
the cardinality of $\{i\}$ is the multiplicity $m_i$, for $i=1,\dots,g$.
By
Proposition~\ref{BC1}, for any point $p\in
M$, the $i^{\text{th}}$ curvature surface through $p$ is either an
open subset of a sphere of dimension $m_i$ in $\R^n$ or an open
subset of a plane of dimension $m_i$.  Inversion in a sphere whose
center does not lie on the union of all these spheres and planes will
transform all curvature surfaces of 
$f(M)$ to open subsets of spheres.  Thus, we may assume that all
curvature surfaces of 
$f$ are open subsets of spheres of the appropriate dimensions.

Denote
the center of the sphere of dimension $m_i$ containing the
$i^{\text{th}}$ curvature surface through $f(p)$ by 
\begin{equation}\label{A8}
c_i(p) \in \R^n,
\end{equation}
and its radius by
\begin{equation}\label{A9}
r_i(p) = |f(p) - c_i(p)|.
\end{equation}
By Proposition~\ref{BC1} and Remark~\ref{AR1}, $c_i(p)$ is determined by
$j^3_p(f)$.  If $e_a(p)$, for $a \in \{i\}$, form an orthonormal
basis 
of principal vectors for the $i^{\text{th}}$ principal curvature, then
they are determined by $j^2_p(f)$.

\begin{lemma}\label{AL2}
For any point $p\in M$, there is an analytic algebraic
paramet\-ri\-zation of an open subset about $f(p)$ of
the $i^{\text{th}}$ curvature surface through
$f(p)$ given by
\begin{equation}\label{A10}
\x(s_i) = c_i(p) + s^0_i(s_i)(f(p) - c_i(p)) + r_i(p) \sum_{a
\in \{i\}} s^a_i e_a(p),
\end{equation}
where 
\begin{equation}\label{A10a}
s^0_i(s_i) = \sqrt{1 - |s_i|^2},
\end{equation}
for all $s_i$ in the 
open ball $B^{m_i}(\epsilon_i)\subset \R^{m_i}$, for some number
$\epsilon_i$ satisfying $0<\epsilon_i<1$.  The components of
$s_i$ are denoted $s_i^a$, for each $a \in \{i\}$, so
\[
|s_i|^2 = 1- \sum_{a\in \{i\}}(s^a_i)^2.
\]
The vectors
$c_i(p)$, $f(p)$, and $e_a(p)$, for all $a \in \{i\}$, are
determined by $j^3_p(f)$.  
\end{lemma}
The proof is contained in Example~\ref{Ex2.1}, which shows that the
map $\x(s_i)$ 
in~(\ref{A10}) is algebraic analytic.  This lemma is essentially
showing that a sphere of any radius, dimension, and center in $\R^n$
is a Nash submanifold (see \cite[Definition 2.9.9, p.\  57]{RAG}).

Since the map $\x(s_i)$ in~(\ref{A10}) is an embedding, there exists
an embedding $x_i:B^{m_i}(\epsilon_i) \to M$ such that $f(x_i(s_i))
= \x(s_i)$, for all $s_i \in B^{m_i}(\epsilon_i)$.

\begin{lemma}\label{AL3} For any $s_i \in B^{m_i}(\epsilon_i)$ in the
parametrization~(\ref{A10}), $j^k_{x_i(s_i)}(f)$ depends 
analytically algebraically on $s_i$ and on $j^{k+3}_p(f)$, for any $k
\geq0$. 
\end{lemma}

\begin{proof} For some $\delta>0$, there exists $y:B^{n-1-m_i}(\delta)
\to Y \subset M$, an embedded 
submanifold through $y(0) = p$ transverse to the $i^{\text{th}}$
curvature surface through each of its points.  Let $c_i(y)$ be the
center of the 
$i^{\text{th}}$ curvature surface through $f(y)$.  Let $e_j$, $j =
1,\dots,n-1$, be a smooth orthonormal frame field of principal
vectors on a neighborhood of $p$.  For each $y \in Y$,
\begin{equation}\label{A13}
f(x_i(s_i),y) = c_i(y) + s^0_i(s_i)(f(y)-c_i(y)) +
r_i(y)\sum_{a\in \{i\}} s_i^a e_a(y)
\end{equation}
is a parametrization of a neighborhood of the $i^{\text{th}}$
curvature surface through $f(y)$.  Thus, for $\epsilon_i$ and $\delta$
sufficiently small, $(x_i,y)$ is a
parametrization of a neighborhood of $p$ in $M$.

From~(\ref{A13}) we see that
the partial derivatives of $f$ with respect to the $s_i$
variables at $(s_i,y)$ depend on
$j^3_{(0,y)}(f)$ and analytically algebraically on $s_i$, since all
partial derivatives of $s_i^0(s_i)$ are analytic algebraic (see
Example~\ref{Ex2.1}).  
Again from~(\ref{A13}) we see that the partial
derivatives of $f$ with respect to the $y$ variables at
$(s_i,y)$ depend analytically algebraically on $s_i$ and on
$j^1_{(0,y)}(c_i)$, 
$j^1_{(0,y)}(f)$, $j^1_{(0,y)}(r_i)$, and $j^1_{(0,y)}(e_a)$, for all
$a \in \{i\}$.  Since $c_i(0,y)$, $f(0,y)$, $r_i(0,y)$, and the
$e_a(0,y)$ are determined by $j^3_{(0,y)}(f)$, it follows that their
1-jets are determined by
\begin{equation}\label{A14}
j^1_{(0,y)}(j^3(f)) = j^4_{(0,y)}(f).
\end{equation}
Taking higher derivatives of~(\ref{A13}) in
this 
way, we see that 
$j^k_{(x_i(s_i),p)}(f)$ depends analytically algebraically on $s_i$ and on
\begin{equation}\label{A15}
j^k_p(j^3(f)) = j^{k+3}_p(f).
\end{equation}
\end{proof}

\begin{theorem}\label{AT2} For any point $p\in M$, there exists an
analytic algebraic paramet\-ri\-zation of a neighborhood of $f(p)$ in
$f(M)\subset 
\R^n$, where analytic algebraic means that the component functions are
analytic algebraic.  This collection of analytic algebraic
parametrizations defines an 
analytic structure on $M$ with respect to which $f$ is analytic.  In
addition, it shows that any point $p\in M$ has an open neighborhood
$U\subset M$ such that $f(U)$ is a semi-algebraic subset of $\R^n$.
\end{theorem}

\begin{proof} The idea is to iterate parametrizations along curvature
surfaces, starting with the first.  Choose a point $p\in M$.  Setting
$i=1$ in equation~(\ref{A10}) gives an analytic algebraic parametrization of
a neighborhood of $f(p)$ of the 
first curvature surface through $p$.  To simplify the following
notation, we will write $f(s_1)$ in place of $f(x_1(s_1))$, and
likewise $c_2(s_1))$, $r_2(s_1))$, and $e_b(s_1)$ without showing
explicitly the composition with the map $x_1(s_1)$, etc.

For each $s_1$, parametrize
a neighborhood of the second curvature surface through $f(s_1)$ by
\begin{equation}\label{A16}
f(s_1,s_2) = c_2(s_1) +
s_2^0(s_2)(f(s_1)-c_2(s_1)) + r_2(s_1) \sum_{b\in
\{2\}}s_2^b e_b(s_1)
\end{equation}
where $c_2(s_1)$, $f(s_1)$, and $e_b(s_1)$, $b \in
\{2\}$, are determined by $j^3_{s_1}(f)$, which in turn depends
analytically algebraically on $s_1$ and $j^3_p(j^3(f)) = j^6_p(f)$, by
Lemma~\ref{AL3}. 

If $g>2$, then one more step should make the
iteration clear.  For each $s_1,s_2 \in B_1(\epsilon_1) \times
B_2(\epsilon_2)$, parametrize a neighborhood of the third 
curvature surface through $f(s_1,s_2)$ by
\begin{equation}\label{A17}
\begin{array}{rl}
f(s_1,s_2,s_3) &= c_3(s_1,s_2) +
s^0_3(s_3) (f(s_1,s_2) -
c_3(s_1,s_2)) \\
&+ r_3(s_1,s_2) \sum_{c \in \{3\}}
s_3^c e_c(s_1,s_2).
\end{array}
\end{equation}
Now $c_3(s_1,s_2)$, $f(s_1,s_2)$, and
$e_c(s_1,s_2)$ are determined by
$j^3_{(s_1,s_2)}(f)$, which in turn depends analytically algebraically
on
$s_2$ and 
\begin{displaymath}
j^3_{(s_1,0)}(j^3(f)) = j^6_{(s_1,0)}(f),
\end{displaymath}
which in turn depends analytically algebraically on $s_1$ and 
\begin{displaymath}
j^3_p(j^6(f)) = j^9_p(f).
\end{displaymath}  

Continuing in this way, we parametrize a neighborhood of the
$g^{\text{th}}$ 
curvature surface through $f(s_1,\dots,s_{g-1})$ by
\begin{equation}\label{A18}
\begin{array}{rl}
f(s_1,\dots,s_g) &= c_g(s_1,\dots,s_{g-1}) \\
&+
s^0_g(s_g)
(f(s_1,\dots,s_{g-1})-c_g(s_1,\dots,s_{g-1})) \\
&+ r_g(s_1,\dots,s_{g-1}) \sum_{d \in \{g\}} s^d_g
e_d(s_1,\dots,s_{g-1}),
\end{array}
\end{equation}
which is analytic algebraic in $(s_1,\dots,s_g) \in
B^{m_1}(\epsilon_1) \times \dots 
\times B^{m_g}(\epsilon_g)$ and depends on the finite set of constants
determined by
$j^{3g}_p(f)$.
By a standard argument (see, for example, Thorbergsson \cite[p.\ 497]{Th1}),
$f(s_1,\dots,s_g)$ parametrizes a neighborhood of $p \in M$.

Finally, since a finite product of open balls is a semi-algebraic
subset of $\R^{n-1}$, its image under the analytic
algebraic map~(\ref{A18}) is
a semi-algebraic subset of $\R^n$.
\end{proof}

\begin{corollary}\label{Main} A connected proper Dupin hypersurface
$M$ in $\R^{n}$
is an open subset of a connected component of 
the irreducible algebraic set $\overline{M}^{\text{zar}}$ of dimension $n-1$.
\end{corollary}
\begin{proof} By the theorem, for any point $p\in M$, there
is an open  
neighborhood $U\subset M$ of the point that is a semi-algebraic subset
of $\R^n$. The result now follows from Lemma~\ref{lemma0}.
\end{proof}

A slightly more general result holds that we shall apply to Dupin
hypersurfaces. 

\begin{corollary}\label{foundation} If a connected analytic
hypersurface $N\subset \R^{n}$ 
contains a connected proper Dupin hypersurface $M$, then $N$ is 
an open subset of the connected component of an
irreducible algebraic set that contains $M$.
\end{corollary}
\begin{proof} An open subset of $M$ is open in $N$, so the result
follows from Lemma~\ref{lemma0}.
\end{proof}

\section{Kuiper's Conjecture}\label{sec:5}
As observed in the Introduction, Kuiper's Conjecture would be proved
if we could prove that any connected non-proper Dupin hypersurface in $\R^n$ is algebraic.

A connected Dupin hypersurface $M$ in $\R^n$ has an open dense subset
${\mathcal G}$ such that each connected  
component of ${\mathcal G}$ is proper Dupin. In fact, ${\mathcal G}$ is the set on which
the multiplicities of the principal curvatures are locally constant (\cite{Singley}).
Decompose ${\mathcal G}$ into its 
at most countably many disjoint connected components
${\mathcal G}_{1}$, ${\mathcal G}_{2}$, ${\mathcal G}_{3},\ldots$.   By Corollary~\ref{foundation}, if
the Dupin hypersurface $M$ were an analytic submanifold of $\R^n$,
then, since it contains the connected proper Dupin hypersurface ${\mathcal G}_1$,
it would be contained in an irreducible algebraic set of dimension $n-1$.
Hence, our goal is to prove that $M$ is an analytic hypersurface of
$\R^n$. 

Let ${\mathcal G}^c$ denote the complement of ${\mathcal G}$ in $M$.  It is the set of
points in $M$ that do not have a proper Dupin neighborhood.  For
example, let $T^2 \subset \R^3$ be a torus of revolution, and let $M^3
\subset \R^4 = \R^3 \times \R$ be a tube over $T^2$ in $\R^4$ of
sufficiently small radius $\epsilon$ so that $M^3$ is an embedded hypersurface in
$\R^4$.  Then $M^3$ is a Dupin hypersurface that is not proper Dupin.
$M^3$ has three distinct principal curvatures at all points except for
the two tori $T^2 \times \{\pm\epsilon\}$ on which there are only two
distinct principal curvatures.  In this case, the set ${\mathcal G}$ is all of
$M^3$ except $T^2\times \{\pm\epsilon\}$, and ${\mathcal G}^c$ equals $T^2 \times
\{\pm\epsilon\}$.  See \cite[p.\ 188]{CecilRyan}.


\begin{definition}\label{finite} A connected Dupin hypersurface $M$ of
$\R^n$ has
the \textit{local finiteness property} if there is a subset $S
\subset {\mathcal G}^{c}$, closed in $M$, such that $S$ disconnects
$M$ into only a finite number of connected components, and for each point
$p\in {\mathcal G}^{c}$ not in $S$, there is an open neighborhood $W$ of $p$
in $M$
such that $W\cap {\mathcal G}$ contains a finite number of connected
open sets whose union is dense in $W$.
We call $S$ the set
of \textit{bad points} in ${\mathcal G}^{c}$.
\end{definition}

For instance, if ${\mathcal G}^c$ is a finite union of compact, connected submanifolds, or
more generally, is a locally finite CW complex, then $M$ has the local
finiteness property.  An example with nonempty $S$ is a set ${\mathcal G}^c$ that
consists of the boundaries of an infinite nested sequence $T_1$,
$T_2, \dots$, of open neighborhoods, with the closure of $T_{n+1}$
properly contained in $T_n$, such that the intersection $S =
\cap_nT_n$ is a submanifold of codimension greater than one in $M$.

\begin{theorem}\label{mainth} Let $M$ be a connected Dupin
hypersurface in $\R^n$.  If $M$ has the local finiteness property, 
then it is an analytic submanifold of $\R^n$, and
is therefore contained in a connected component of dimension $n-1$ of
an irreducible algebraic set. 
\end{theorem}

\begin{proof}   Let ${\mathcal G}\subset M$ be the open dense subset of $M$
whose connected components are proper Dupin hypersurfaces in $\R^n$.
Let $S \subset {\mathcal G}^c$ be the set of bad points in ${\mathcal G}^c$.
By the local finiteness property, each point $p\in {\mathcal G}^{c}$ not in $S$
has a neighborhood $W$, open in $M$, such that $W \cap {\mathcal G}$ contains a finite number of connected open sets
$U_{1},\ldots,U_{s}$, whose union is dense in $W$. 
Note that ${\mathcal G}$ dense in $M$ implies that
\[
W \subset \cup_1^s \overline{U}_i,
\]
where $\overline{U}_i$ is the topological closure of $U_i$ in $M$.
By Corollary~\ref{foundation}, each $U_{i}$ is contained in  
an irreducible algebraic set $C_{i}$.
Then 
$\overline{U}_{i}$ is a subset of $C_{i}$ since $C_{i}$ is closed.
Hence, 
\[
W \subset \cup_1^s \overline{U}_i \subset \cup_1^s C_i,
\]
which is a semi-algebraic subset of $\R^n$.

By Lemma~\ref{keylemma},  since $W$ is a $C^\infty$-manifold,
there is a connected open semi-algebraic subset $U$ of
$\cup_{1}^{s} C_{i}$ contained in $W$ with $p\in U$.  
The intersection $B$ of
$U$ with any open ball of $\R^n$ centered at $p$ is still an open
subset of $M$ and a semi-algebraic subset.  Thus,
we may assume $B$ is so small that it is the graph of 
$h:D\rightarrow\R$ for some open set $D\subset \R^{n-1}$
by performing a linear change of coordinates $(x^{1},\ldots,x^{n})$ 
in $\R^{n}$ if necessary.
Since the projection
$$
\pi:(x^1,\ldots,x^{n})\mapsto (x^{1},\ldots,x^{n-1}),
$$ 
which is semi-algebraic, sends 
$B$ to $D$, 
we see $D$ is a semi-algebraic subset of $\R^{n-1}$.
Thus, the map $h$ 
is semi-algebraic, because its graph 
$$
\{(x,h(x)): x\in D\}=B
$$
is a semi-algebraic subset of $\R^{n}$.
Therefore, $h$ is a Nash function and so it must be
analytic algebraic as noted in Proposition~\ref{char}, and the map
\[
H:D \to B, \; H(x) = (x,h(x))
\]
is an analytic algebraic parametrization of the open neighborhood $B
\subset M$ about $p$.  Since $p\in {\mathcal G}^c\setminus S$ was arbitrary, and
since every 
point of ${\mathcal G}$ has a neighborhood with 
an analytic algebraic parametrization by
Theorem~\ref{AT2}, it follows that $M\setminus S$ is an analytic
submanifold of $\R^n$. 
Since $M\setminus S$ has a finite number of connected components, each of
which contains an open subset that is a semi-algebraic subset of $\R^n$,
it follows that $M\setminus S$ is contained in the union $C$ of finitely
many irreducible algebraic sets.
Then $M$ is contained in the topological closure of
$M\setminus S$ which in turn must be contained in the closed set $C$.

Now let $q$ be any point in $S$ and consider $M$ itself as an open
neighborhood about $q$.  Then $M$ is
a $C^\infty$ manifold contained 
in the semi-algebraic set $C$.  Apply the above argument to conclude
that some open neighborhood of $q$ in $M$ is the image of an analytic
algebraic parametrization.  Hence, every point of $M$ has an analytic
parametrization, and $M$ is an analytic hypersurface of $\R^n$.
\end{proof}

\begin{question}\label{question} Is the local finiteness property true
for any 
connected Dupin hypersurface in $\R^n$?
\end{question}

If the answer is yes, then we could resolve Kuiper's conjecture for
taut submanifolds of $\R^n$ in the following manner.

First note that Pinkall \cite{Pinkall2} extended the notion of Dupin to submanifolds in $\R^n$ of codimension
greater than one as follows.  Let $M$ be a submanifold in $\R^n$ of codimension greater than one, and let
$UN(M)$ denote the unit normal bundle of $M$ in $\R^n$.  A {\em curvature surface} of $M$ is a connected submanifold 
$S \subset M$ for which there is a parallel (with respect to the normal connection) section $\eta:S \rightarrow UN(M)$ 
such that for each $x \in S$, 
the tangent space $T_xS$ is equal to some  principal space of the shape operator $A_{\eta (x)}$.
The corresponding principal curvature
$\kappa :S \rightarrow {\bf R}$ is then a smooth function on $S$.  Pinkall calls a submanifold
$M$ of codimension greater than one {\em Dupin} if along each curvature surface,
the corresponding principal curvature is constant.  This definition is equivalent to Pinkall's \cite{Pinkall}
Lie sphere geometric definition of the Dupin property in terms of the Legendre lift of $M$, and of course,
it agrees with the usual definition of a Dupin hypersurface in the case of codimension one.

Pinkall \cite{Pinkall2} then proved that a taut submanifold of $\R^n$ must be Dupin, although not necessarily proper
Dupin.  He also proved that if
$M$ is an embedded taut compact submanifold in $\R^n$ and $M_\epsilon$ is a tube over $M$ of sufficiently small
radius that $M_\epsilon$ is an embedded hypersurface in $\R^n$, then $M$ is Dupin if and only if $M_\epsilon$ is Dupin
(since their Legendre lifts are Lie equivalent \cite{Pinkall}), and $M$ is taut (with respect to $\Z_2$-homology)
if and only if $M_\epsilon$ is taut.

We now prove that a taut submanifold $M$ in $\R^n$ is algebraic if the local finiteness property holds on the tube
$M_\epsilon$.

\begin{theorem}\label{Kuiper} 
{\rm (a)} If $M^{n-1}$ is a connected, compact taut hypersurface in $\R^n$ that satisfies the local finiteness property,
then $M^{n-1}$ is an analytic submanifold and a connected component of an irreducible algebraic subset of $\R^n$.\\
{\rm (b)} If $M^m$ is a connected, compact taut submanifold of codimension greater than one in
$\R^n$ such that the tube $M_\epsilon$ over $M$, which is Dupin,
satisfies the local finiteness property, 
then $M$ is an analytic submanifold and a connected
component of an irreducible
algebraic subset of $\R^n$.
\end{theorem}

\begin{proof}
(a) Since a taut hypersurface must be Dupin, this follows immediately from Theorem~\ref{mainth}.\\
(b) Let $M$ be an embedded taut connected, compact submanifold in $\R^n$ and let $M_\epsilon$ be a tube over $M$ of 
sufficiently small
radius that $M_\epsilon$ is an embedded hypersurface in $\R^n$.  Then $M_\epsilon$ is a Dupin hypersurface by
the work of Pinkall \cite{Pinkall2} mentioned above.  Thus $M_\epsilon \subset \R^n$ is semi-algebraic by
Theorem~\ref{mainth}, since $M_\epsilon$ satisfies the local finiteness property.
Consider the focal map $F_\epsilon : M_\epsilon \rightarrow M \subset \R^n$ given by
\begin{displaymath}
F_\epsilon(x) = x - \epsilon \xi,
\end{displaymath}
where $\xi$ is the outward field of unit normals to the tube $M_\epsilon$.  Any point of $M_\epsilon$ has an open
neighborhood $U$ parametrized by an analytic algebraic map.  By Proposition \ref{good}, the first derivatives
of this parametrization are also analytic algebraic, and thus the Gram-Schmidt
process applied to these first derivatives and some constant non-tangential vector produces the
vector field $\xi$ and shows
that $\xi$ is analytic algebraic on $U$.  Hence $F_\epsilon$ is analytic algebraic on $U$ and so the image 
$F_\epsilon (U) \subset M$ is a semi-algebraic subset of $\R^n$.  Covering $M_\epsilon$ by finitely many sets
of this form $U$, we see that $M$, being the union of their images under $F_\epsilon$, is a semi-algebraic subset
of $\R^n$. Then $\overline{M}^{\text{zar}}$ is an irreducible algebraic of the same dimension as $M$ and contains 
$M$.  
\end{proof}

The preceding proof uses the finiteness condition on the tube $M_\epsilon$ to prove that $M_\epsilon$
is analytic via Theorem~\ref{mainth}.  This will also follow from the assumption that $M$ itself is analytic.
Thus we have the following result.

\begin{corollary}\label{Cor:5:1} If $M$ is an analytic connected, compact taut submanifold of
$\R^n$, 
then $M$ is a connected
component of an irreducible
algebraic subset of $\R^n$.
\end{corollary}

\section{Alexander Cohomology}\label{sec:6}
We now begin our work to show that the local finiteness property is satisfied
by taut hypersurfaces satisfying certain additional restrictions on their dimensions
and the multiplicities of their principal curvatures.  A principal tool in our work is the
theory of Alexander cohomology for noncompact manifolds.  We begin with the following
definition.

\begin{definition} 
\label{end1}
An \textit{end} of a noncompact manifold $X$ is an equivalence class of sequences of
connected open neighborhoods $X, U_1, U_2,\ldots$, where $U_{n+1}\subset
U_{n}$, such that the intersection of the closures of these sets is the empty set,
subject to the equivalence relation $(X,U_1,U_2,\ldots)\sim
(X,W_1,W_2,\ldots)$ if for each $U_{i}$ there is a $j$ with
$U_{i}\subset W_{j}$, and for each $W_{k}$ there is an $m$ such that
$W_{k}\subset U_{m}$.
\end{definition}

Though the definition is intuitive, another equivalent intrinsic definition is often
more convenient to work with.

\begin{definition}
\label{end2} 
An \textit{end} of a noncompact manifold $X$ is a function $\epsilon$ which assigns each compact
set $K\subset X$ to a subset of $X$ such that $\epsilon(K)$ is a connected component of
$X\setminus K$ for each $K$, and if $K\subset L$ then $\epsilon(L)\subset\epsilon(K)$.
\end{definition}

For a proof that these two definitions are equivalent, see \cite[Proposition 9, p.\ 9]{HR}.

We can ``glue'' the ends of $X$ to $X$ to make the resulting space compact (see \cite{Ray}).
In particular, the set of ends is not empty by the noncompactness of $X$.

\begin{definition} Let $S$ be a closed subset of a manifold $M$.
The $k$-th \textit{Alexander cohomology group} is defined to be
$$
\overline{H}^{k}(S):=\lim_{\rightarrow} H^{k}(V),
$$
where $H^{k}(V)$ is the
usual $k$-th singular cohomology group of $V$ (with, say, the coefficient
ring $\Z_2$; in fact, all commutative rings with identity suffice), and
$\lim_{\rightarrow}$ denotes taking the direct limit over the directed
system of open sets $V$ containing $S$ in the ambient space $M$.
\end{definition}

In other words, we define $\overline{H}^{k}(S)$ to be the set of equivalence
classes of the disjoint union of $H^{k}(V)$, for all open $V$ containing $S$,
where $a_{U}\in H^{k}(U)$ and $a_{V}\in H^{k}(V)$ are considered equivalent
if there is a third open set $W$ containing $S$, where $W$ is contained in both
$U$ and $V$, with the inclusion maps $\iota^{U}_{W}:W\rightarrow U$ and
$\iota^{V}_{W}:W\rightarrow V$, and the pullback morphisms
$(\iota^{U}_{W})^{*}: H^{k}(U)\rightarrow H^{k}(W)$ and 
$(\iota^{V}_{W})^{*}: H^{k}(V)\rightarrow H^{k}(W)$, such that
$(\iota^{U}_{W})^{*}(a_{U})=(\iota^{V}_{W})^{*}(a_{V})$.

Let us look at $\overline{H}^{0}(S)$ to gain some insight.
First observe that if $W\subset U$ are two open neighborhoods of
$S$, then the 
morphism $(r^{U}_{W})^{*}:H^{0}(U)\rightarrow H^{0}(W)$ induced by
the inclusion $r^{U}_{W}: W\rightarrow U$ is nothing other than
restricting an element $a_{U}\in H^{0}(U)$, which is a $\Z_2$-valued
locally constant function, to the same function on $W$. Hence, $a_{U}\in H^{0}(U)$ and
$a_{V}\in H^{0}(V)$ are equivalent in the definition of the Alexander
cohomology group $\overline{H}^{0}(S)$ if and only if there is a smaller open
neighborhood $W$ of $S$ contained in both $U$ and $V$ such that the two
$\Z_2$-valued functions
$a_{U}$ and $a_{V}$ restrict to the same function on $W$. In other words,
$\overline{H}^{0}(S)$ is the ring of germs of locally constant $\Z_2$-valued functions
around $S$. 

Here, for any topological space $X$, a function $f$ from $X$
to $\Z_2$ is said to be {\em locally constant} if there is an open
covering ${\mathcal U}$ of $X$
such that $f$ is constant on each element of ${\mathcal U}$.

\begin{theorem}\label{alex}
$\overline{H}^{0}(S)$ is the
$\Z_2$-module of locally constant functions from $S$ to $\Z_2$.
In particular, if all the (topological) connected components of $S$ are open in $S$,
then the rank of $\overline{H}^{0}(S)$ is identified with
the cardinality of the (topological) connected components.
\end{theorem}

For a proof, see Theorem 5 and Corollary 6 in \cite[p.\ 309-310]{Span}.

\begin{corollary}
\label{corollary} 
Let $M$ be a compact and simply connected
manifold. Let $S\subset M$ be a compact subset which
disconnects $M\setminus S$ into $\alpha$ connected components. Suppose all the
connected components of $S$
are open (whose cardinality is necessarily a finite number $m$) in $S$.
Then the number of ends $\epsilon$ of $X=M\setminus S$ is finite, and 
\begin{equation}\label{count}
\epsilon=\alpha+m-1.
\end{equation}
\end{corollary}
\begin{proof} 
First of all, observe that the number of (topological)
connected components of $S$ 
is finite, since $S$ is the disjoint union of these open components and $S$
is compact; let this number be $m$.
Therefore, by Theorem~\ref{alex}, the rank of $\overline{H}^{0}(S)$ is $m$.

We claim that the number $\epsilon$ of ends of $X$ is finite.
Suppose $\epsilon=\infty$. Then
there is a decreasing sequence of open sets $W_{1},W_{2},\ldots$ in $X$
such that $X \setminus W_{j}$ is compact, $\cap_{j} \overline{W}_{j}=\emptyset$, and the number of unbounded
connected components of $W_{j}$ is an increasing sequence
diverging to infinity ($\overline{W}$ is the closure of $W$).
(By \cite[Theorem 3.9, p.\ 111]{HY}, each of $W_{j}$
has finitely many unbounded connected components, since 
$X\setminus W_{j}$ is compact.)
Let $n_{j}$ be the number of connected components of $W_{j}$. Let
$U_{j}=W_{j}\cup S$. Then
$U_{j}$ is a neighborhood of $S$ in $M$ since $U_{j}=M\setminus (X \setminus W_{j})$.
Let $m_{j}$
be the number of connected components of $U_{j}$ in $M$.

Observe that since $M$ is covered by the open sets $M \setminus S$ and $U_{j}$, the
Mayer-Vietoris
sequence, with $H^{1}(M)=0$ by simple connectedness, says
$$
0\leftarrow H^{0}(U_{j}\setminus S)\leftarrow
H^{0}(M\setminus S)\oplus
H^{0}(U_{j})\leftarrow H^{0}(M)\leftarrow 0
$$
so that
\begin{eqnarray}\label{estimate}
n_{j}&=&\dim H^{0}(U_{j}\setminus S)=\dim H^{0}(M\setminus S)+\dim H^{0}(U_{j})-1\\
&=&\alpha+m_{j}-1.\nonumber
\end{eqnarray}

It follows that the sequence $\{m_j\}$ also diverges to $\infty$. Now the direct limit
in the definition of Alexander cohomology can be taken over the directed
system $U_{1},U_{2},\ldots$ in place of the system of all open neighborhoods of $S$.
This is because for each open neighborhood $U$ of $S$, there is a
$W_{j}\subset U \setminus S$, and hence a $U_{j}\subset U$, such that if we define
$a_{_{j}}\in H^{0}(U_{j})$ to be the
restriction of
$a_{U}\in H^{0}(U)$ on $U_{j}$, then $[a_{U}]=[a_{U_{j}}]\in \overline{H}^{0}(S)$.

The restriction maps $r_{j}:U_{j+1}\rightarrow U_{j}$ induce the
restriction morphisms
$$
(r_{j})^{*}:H^{0}(U_{j}) \rightarrow H^{0}(U_{j+1}).
$$
Moreover, $a_{j}\in H^{0}(U_{j})$ is identified with $a_{j+1}\in H^{0}(U_{j+1})$
if and only if $a_{j+1}=(r_{j})^{*}(a_{j})$.
Then
$$
(r_{j})^{*}:\oplus \Z_2 \rightarrow \oplus \Z_2,
$$
where the domain ring consists of $m_{j}$ copies of $\Z_2$ and
the target ring of $m_{j+1}$ copies with $m_{j}<m_{j+1}$.

Starting with $H^{0}(U_{1})$, which has cardinality $2^{m_{1}}$, $(r_{1})^{*}(U_{1})$
is of cardinality no greater than $2^{m_{1}}$ in the space $H^{0}(U_{2})$ of
cardinality $2^{m_{2}}>2^{m_{1}}$,
$(r_{2})^{*}(r_{1})^{*}(U_{1})$ is of cardinality at most $2^{m_{1}}$ in the space $H^{0}(U_{3})$
of cardinality $2^{m_{3}}>2^{m_{1}}$, etc. Now pick one element $a_{2}$ in $H^{0}(U_{2})$
not in
the image of $H^{0}(U_{1})$ under $(r_{1})^{*}$, an element $a_{3}$ not in the
image of $H^{0}(U_{1})$ under $(r_{2})^{*}(r_{1})^{*}$, etc. Then
$(a_{2},a_{3},\cdots)$ is not equal to any element in $\overline{H}^{0}(S)$
whose representative comes from $H^{0}(U_{1})$. Continuing in this fashion,
we can construct an infinite number of distinct elements in
$\overline{H}^{0}(S)$,
since $m_1<m_2 <\cdots$ increasingly diverge to infinity.  The upshot is that
this implies that the cardinality of $\overline{H}^{0}(S)$ is infinite, which contradicts the fact
that it is equal to $2^{m}$. 

The contradiction proves that the number $\epsilon$ of ends is finite.
It follows that, since $n_{j}$ is the number of connected components of $W_{j}$, then
$n_{j}=\epsilon$ for sufficiently large $j$, so
that~\eqref{estimate} gives
$$
m_{j}=1+\epsilon-\alpha
$$
for sufficiently large $j$.

In the case when $m=1$, i.e., when $S$ is connected, we have
$m_{j}=1$ for large $j$, since
the closure of each of the finite ends intersects $S$ (or else it would not be an end), so that in fact
the union of each end with $S$ is connected, and so $U_{j}$, being the union of all $W_{j}$
and $S$, with $W_{j}$ composed of the ends, is connected. 
In general, since all the components of $S$ are open, we can
separate them by disjoint open sets in $M$ and hence can work on each
component to
conclude that in fact $m_{j}=m$ for large $j$, from which we derive
$$
\epsilon=\alpha+m-1.
$$
\end{proof}
We assume in the preceding corollary that $H^{1}(M)=0$ only for simplicity.
When $H^{1}(M)\neq 0$, we have the connecting homomorphism
$$
\leftarrow H^{1}(M)\leftarrow H^{0}(U_{j}\setminus S)\leftarrow
$$
when chasing upward the Mayer-Vietoris sequence, so that in fact
\begin{equation}\label{stability}
n_{j}=\alpha+m_{j}-1+\delta_{j}
\end{equation}
with $\delta_{j}$ the dimension of the image of $H^{0}(U_{j}\setminus S)$ via the
connecting homomorphism.
Clearly
$$
\delta_{j}\leq b_{1},
$$ where $b_{1}$ is the first Betti number of $M$.

By the same reasoning as above we see that $n_{j},m_{j},\delta_{j}$ eventually
stabilize, so that we have
\begin{equation}\label{topology}
\epsilon=\alpha+m-1+\delta
\end{equation}
with $\delta\leq b_{1}$. We conclude the following.

\begin{corollary}\label{bound} Let $S$ be a closed set in
a compact
manifold $M$. Suppose all the (topological) connected components of $S$
are open (and so necessarily of finite cardinality $m$) in $S$, 
and suppose the number of ends of $X=M\setminus S$ is finite with cardinality
$\epsilon$.
Then the number of connected components of $X$ is
no greater than $\epsilon-m+1$.
\end{corollary}

\begin{proof} In Equation~\eqref{topology}
$$
\alpha=\epsilon-m+1-\delta\leq\epsilon-m+1.
$$
\end{proof}

\begin{remark}\label{fractal} It is important to observe that we calculated Alexander cohomology
under the assumption that the connected components of $S$ are all open in $S$. Conversely,
when we know in advance that the number of ends of $X$ is finite, we
can assert that the (topological) connected components of $S$ are all open
in $S$. (As long as $S$ is closed, $M$ need not be compact.) This is
because~\eqref{stability}
and the finite cardinality of ends imply that the numbers $m_{j}$
are always bounded above by a
constant $\rho$. (Note, in particular, that the number of connected
components of $X$, which is $\alpha$ in~\eqref{stability}, is then always finite.) Fix a large $j$ in the process and 
suppose the open set $U_{j}$
decomposes into $\beta$ connected components $U_{j1},\cdots,U_{j\beta}$. We may
assume all $S_{k}:=U_{jk}\cap S\neq\emptyset,1\leq k\leq\beta$; or else we
ignore those having empty intersection with $S$ and redefine $U_{j}$ to be the union
of those intersecting $S$. Then for $U_{j+1}\subset U_{j}$, since
$S\subset U_{j+1}$, none of the $U_{j+1,k}:=U_{j+1}\cap U_{jk},1\leq k\leq\beta,$
are empty. It follows that $U_{j+1}$ consists of at least $\beta$ connected components.
Hence $m_{j}\leq m_{j+1}\leq \cdots\leq\rho$. This implies that eventually all large
$m_{t}$ are equal to a fixed number $m$. Fix such $t$.
The $m$ connected components of $U_{t}$ intersect $S$ in $m$ nonempty
disjoint 
sets
$S_{1},\ldots,S_{m}$, which in turn are contained in the $m$
connected components of all $U_{r},r\geq t,$ by the process. Note that the restriction morphisms $(r_{t})^{*}$
are all isomorphisms for large $t$. Hence the rank of $\overline{H}^{0}(S)$
is $m$, which is the rank of the module of locally constant functions on $S$.
Now the sets $S_{1},\ldots,S_{m}$ are both open and closed in $S$, and they form an open covering.
It follows that the locally constant functions on $S_{1},\ldots,S_{m}$ already
exhaust $\overline{H}^{0}(S)$. Each $S_{j}$ must be connected, or else, we could cut it up into two disjoint nonempty
open sets to find a refined open covering with $m+1$ open sets, contradicting the
rank of $\overline{H}^{0}(S)$ being $m$. Hence all $S_{1},\ldots,S_{m}$ are open 
connected components.
\end{remark}

\begin{definition}
The set $S$ is called the \textit{end set} of $X=M\setminus S$. If $X$ has $h$ ends, then there is
an open neighborhood of $U$ of $S$ such that for any open neighborhood
$V$ of $S$ contained in $U$, $V\setminus S$ has $h$ (unbounded) connected components,
each of
which corresponds uniquely to an end. We call these components the \textit{end components} of $X$.
\end{definition}

\section{Application to Compact Taut Submanifolds}\label{sec:7}

We begin this section by recalling some fundamental facts about taut embeddings.  
An embedded compact, connected submanifold $M \subset \R^n$ is
taut if every nondegenerate Euclidean distance function, 
\[
L_p:M \to \R, \quad L_p(x) = d(x,p)^2, \quad p \in \R^n,
\]
has $\beta (M,\Z_2)$ critical points on $M$, where $\beta (M,\Z_2)$ is the sum of the $\Z_2$-Betti numbers
of $M$.  That is, $L_p$ is a perfect Morse function on $M$.  

Kuiper \cite{Kuiper1} showed that a compact, connected submanifold $M \subset \R^n$
is taut if and only if $M$ satisfies the following condition: for every closed ball $B \subset \R^n$ the homorphism
on homology
\begin{equation}
\label{eq:7.1}
H_*(B \cap M, \Z_2) \rightarrow H_*(M,\Z_2)
\end{equation}
induced by the inclusion of $(B \cap M) \subset M$ is injective for all $r \in {\bf R}$.

Tautness is invariant under M\"{o}bius transformations of ${\bf R}^n \cup \{\infty \}$.  Further, a compact submanifold
$M \subset \R^n$ is taut if and only if the embedding $\sigma (M) \subset S^n$, where
$\sigma: \R^n \rightarrow S^n - \{P\}, P \in S^n$, is stereographic projection, is taut in $S^n$, where spherical distance
functions $d_p$ are used instead of Euclidean distance functions $L_p$.  In this section, we will consider submanifolds
of $S^n$ instead of $\R^n$ for various reasons, including the fact that a focal point cannot vanish to infinity in $S^n$.

Since a spherical distance function $d_p (q) = \cos^{-1} (p \cdot q)$ has the same critical points as
the Euclidean height function
$\ell_p (q) = p \cdot q$, for $p,q \in S^n$, a compact submanifold $M \subset S^n$ is taut if and only if it is
{\em tight}, i.e., every nondegenerate height function $\ell_p$ has $\beta (M, \Z_2)$ critical points
on $M$. It is often simpler to use height functions rather than spherical distance functions when studying tautness
for submanifolds of $S^n$, and we will use whichever type of function is most convenient for our argument. 

We now come to a fundamental result on taut submanifolds due to Ozawa \cite{Ozawa}.

\begin{theorem}[Ozawa]
\label{Ozawa}
Let $M$ be a taut compact connected submanifold embedded in $S^n$, and let $\ell_p, p \in S^n$, be a linear height
function on $M$.  Let $x \in M$ be a critical point of $\ell_p$, and let $S$ be the connected component of the
critical set of $\ell_p$ that contains $x$. Then $S$ is \\
{\rm (a)} a smooth compact manifold of dimension equal to the nullity of the Hessian of $\ell_p$ at $x$;\\
{\rm (b)} nondegenerate as a critical manifold;\\
{\rm (c)} taut in $S^n$.
\end{theorem}
We call such a connected component of the critical set of $\ell_p$ a {\em critical submanifold} of $\ell_p$.

\begin{remark}
\label{Terng-Thorbergsson}
Terng and Thorbergsson \cite[p.190]{T-Th} generalized the notion of tautness to submanifolds of arbitrary complete
Riemannian manifolds and proved an analogue of Ozawa's Theorem in that context.
\end{remark}
 
Let $M$ be a compact taut hypersurface in $S^{n}$. 
Consider the
normal exponential map $E:M\times (-\pi,\pi)\rightarrow S^{n}$, where 
$$
E:(p,t)\mapsto \cos(t) p + \sin(t) {\bf n},
$$
with ${\bf n}$ the chosen unit normal field of $M$.
Here, $E$ is smooth at points where $t=-\pi,\pi$.

A point $q = E(p,t)$ is called a {\em focal point of multiplicity} $m>0$ {\em of} $M$ {\em at} $p$ if the nullity of the
derivative $E_*$ is equal to $m$ at $(p,t)$.
The set of all focal points
is the {\em focal set} ${\mathcal F}$ of $M$.
The focal points at $p$
are antipodally symmetric on the circle $E(p,t)$ with each principal curvature
of the form $\cot(t)$ for some $t\neq 0$.
Let ${\mathcal Z}$ be the complement of ${\mathcal F}$ in $S^{n}$.

\begin{lemma}\label{Hausdorff}
${\mathcal Z}$ is connected in $S^{n}$.
\end{lemma}

\begin{proof}
By Federer's version of Sard's theorem \cite[Theorem 3.4.3, p.316]{Fed}, the image of the
critical points of a given smooth function
$f:\R^t\rightarrow \R^s$, at which the rank of the derivative is less than or equal to
$\nu$, is of ${\mathcal H}^{\nu}$-measure 0, where ${\mathcal H}^{\nu}$
denotes the Hausdorff $\nu$-dimensional measure.

Label the principal curvature functions by
$\lambda_{1}\leq\lambda_{2}\leq \cdots\leq\lambda_{n-1}$.
It is known that $\lambda_{j}$
are Lipschitz-continuous, because the principal curvature functions on the linear space ${\mathcal S}$
of all symmetric matrices are Lipschitz-continuous by general matrix
theory (see, for example, \cite[p.\ 64]{Bha})
and the fact that the Hessian of $M$ is a smooth function from $M$ into ${\mathcal S}$.
Let $\cot(t_{j}):=\lambda_{j}, 0<t_{j}<\pi$.
We know that the focal maps
\begin{equation}\label{equation}
f_{j}(p)=\cos(t_{j}) p + \sin(t_{j}) {\bf n}
\end{equation}
are Lipschitz-continuous and are smooth on a dense open set.

For each $i$, let $O_{i}$ be the open subset of $M$ on which $\lambda_{i}$ has
multiplicity $1$. ($O_{i}$ could be empty.) Then $O_{i}$ consists of countably many
open components $O_{ij},j=1,2,\ldots$, such that the restriction of $f_{i}$
to $(O_{ij})$ is an embedded submanifold
of dimension $n-2$. This follows from the fact that $\lambda_{i}$ is constant
on its lines of curvature on each
$O_{ij}$ by the Dupin condition.

Let $Z_{i}$ be the complement of $O_{i}$ in $M$. ($Z_{i}$ could be empty.)
At each point $p$ of $Z_{i}$, the principal curvature $\cot(t_{i})$ must have
multiplicity at least $2$, and so the normal exponential map $E$ has rank
$\leq n-2$ at the point $(p,t_{i})$. Thus the focal point $f_{i}(p)$ lies in the
singular value set $K$ of points for which the derivative of $E$ has rank
$\leq n-2$.

We conclude that the entire focal set ${\mathcal F}$ is composed of the countably many
embedded submanifolds $f_{i}(O_{ij})$ of dimension $n-2$, their antipodal sets
and the set $K$
which has Hausdorff ($n-2$)-measure zero by Federer's theorem quoted above.
Therefore, the Hausdorff ($n-1$)-measure of the whole focal set ${\mathcal F}$
is zero, which implies that the complement ${\mathcal Z}$ of ${\mathcal F}$ is connected 
(see \cite[p.\ 269]{SY}).
\end{proof}

Away from ${\mathcal F}$, the
map $E$ is a local diffeomorphism. For each $0\leq m\leq n-1$, we let
$W_{m}$ be the set of points $(p,t)$ in $E^{-1}({\mathcal Z})$
for which the spherical distance function $d_{q}$, where $q=E(p,t)$, has
index $m$ at $p$. Then the tautness of $M$ implies that
$$
E:W_{m}\rightarrow {\mathcal Z}
$$
is a degree $b_m$ (regular) covering map, where $b_{m}$ is the $m$-th Betti number. Therefore, by the connectedness of
${\mathcal Z}$, the set $W_{m}$ decomposes into
finitely many (connected) covering sheets $W_{m1},W_{m2},\ldots,W_{ms_{m}}$
onto ${\mathcal Z}$, where $s_m \leq b_m$. Note, in particular, that $W_{m}\neq\emptyset$ if and
only if $b_{m}\neq 0$. 
Further, since $W_{0}$ and $W_{n-1}$ are both strips around the $0$-section $M\times \{0\}$,
we will not consider $m=0$ or $m=n-1$. Thus, we will consider only those $m$ with $0 < m < n-1$, such that
$b_m \neq 0$.  For such values of $m$, $W_{m}$
will be composed of connected components either in $M\times (0,\pi)$ or
$M\times (-\pi,0)$. We denote by $W_{mj}^{+}$ those $W_{mj}$ contained in
$M\times (0,\pi)$. We set 
$W_{m}^{+}:=\cup_{j} W_{mj}^{+}$. 

\begin{lemma}\label{importantlemma} 
All the (topological) connected components of ${\mathcal F}$
are open in ${\mathcal F}$. In particular, the number of ends of ${\mathcal Z}$ is finite.
\end{lemma}

\begin{proof}
Since each $f_j(M),1\leq j\leq n-1,$ defined in~\eqref{equation} is compact
and path-connected
in $S^{n}$, we can group them and the sets antipodal to them into classes where the union 
of those $f_j(M)$ in each class is path-connected whereas the unions from different classes
are disjoint. Call these disjoint unions $X_{1},\ldots,X_{r}$, each of which,
being a finite union of compact sets,
is closed in $S^{n}$ and path-connected. Hence by the Urysohn separation
lemma, there are disjoint open sets $O_{1},\ldots,O_{r}$ that
contain $X_{1},\ldots,X_{r}$, respectively. This means that each (topological)
component, now being just the path-connected components, are open in the relative topology. The conclusion
follows by Corollary~\ref{corollary} since the ambient sphere is simply
connected and ${\mathcal F}$ does not disconnect $S^{n}$ so that $\alpha=1$
in~\eqref{count}. The number of ends of ${\mathcal Z}$ equals the number of path-connected components
of ${\mathcal F}$.
\end{proof}

\begin{corollary} Each $W_{mj}^{+}$ has a finite number of ends.
\end{corollary}

\begin{proof} This follows since $E:W_{mj}^{+}\rightarrow {\mathcal Z}$
is a proper map because it is a covering map of finite degree and ${\mathcal Z}$ has a
finite number of ends.
\end{proof}

Let $pr:M\times (-\pi,\pi)\rightarrow M$ be the projection. Since $pr$
is an open map, the sets
\begin{equation}\label{pr-eq}
U_{mj}^{+}:=pr(W_{mj}^{+})
\end{equation}
are open and connected in $M$. We also set
\begin{equation}\label{union}
U_{m}^{+}:=\cup_{j} U_{mj}^{+} = pr(W_{m}^{+}).
\end{equation}
It is clear that $W_{m}^{-}$ and $U_{m}^{-}$ can be similarly defined on $M\times (-\pi,0)$.

\begin{definition}\label{*set} We define $(U_{m}^{*})^{+}$ to
be the collection of all $x$ for which there exists a $t>0$ such that $(x,t)$ is
a regular point of the normal exponential map $E$, and the spherical distance function $d_{y}$, where $y=E(x,t)$,
has index $m$ at $x$.
\end{definition}

Note that $U_{m}^{+}\subset (U_{m}^{*})^{+}$ by definition. Further,
any point $p$ in the complement $(U_{m}^{*})^{+} \setminus U_{m}^{+}$ satisfies the condition that every
$q=E(p,t),t>0,$ on the normal exponential circle, such that $d_{q}$
has a nondegenerate critical point of index $m$ at $p$, is a focal point of some other point in $M$. 

\begin{lemma} \label{lemma*}
$(U_{m}^{*})^{+}$ is open in $M$.
\end{lemma}

\begin{proof} This follows from the property that
$E$ is a local diffeomorphism around $(x,t),t>0$, so that $E(x,t)$ is
nonfocal along the normal exponential circle.
Hence, for a point $(x',t')$ near $(x,t)$, the point $E(x',t')$ is also
nonfocal along the respective normal exponential circle.
Now,
$d_{y}$ has index $m$ at the nondegenerate critical point $x$. When $x'$ is
sufficiently close to $x$, the function $d_{y'},y'=E(x',t'),$ is a slight perturbation of
$d_{y}$.
Since $x'$ is a nondegenerate critical point of $d_{y'}$ and since 
nondegenerate critical points are locally structurally stable, $d_{y'}$ must
also have index $m$ at $x'$, so that $x'\in (U_{m}^{*})^{+}$. 
\end{proof}

We introduce a space slightly larger than $W_{mj}^{+}$ as follows. Observe
that any two points $(p,t_{1}),(p,t_{2})$ in $W_{mj}$ lie on an interval in the set $\{p\} \times (0,\pi)$ between two
adjacent critical points of the map $E$.
However, it should be noted that, in general, there might exist points $(p,t)$ on the
same interval not belonging to $W_{mj}^{+}$, because $E(p,t)$ is a focal point of some other point in $M$.
In accordance, for
$(p,t_{1})$ and $(p,t_{2})$ in $M\times (0,\pi)$, we say $(p,t_{1})$ is
equivalent to $(p,t_{2})$, denoted $(p,t_{1})\sim (p,t_{2})$, if
the distance functions $d_{q_{1}}$, for $q_{1}=E(p,t_{1})$, and $d_{q_{2}}$, for $q_{2}=E(p,t_{2}),$
have nondegenerate critical points of the same index $m$ at $p$.
We then let $L_{mj}$ be the (trivial) line bundle

\begin{equation}
L_{mj}:=\{(p,t)\in M\times(0,\pi):(p,t)\sim (p,t_{0})\  {\mbox {\rm for some }}(p,t_{0}) \in W_{mj}\}
\end{equation}                                                 
over $U_{mj}^{+}$.

\begin{lemma}\label{finite-ends} $L_{mj}\setminus W_{mj}^{+}$ is of Hausdorff codimension at least $2$
in $L_{mj}$. In particular,
$L_{mj}$ and $U_{mj}^{+}$ are of finite number of ends.
\end{lemma}

\begin{proof} By the openness of $W_{mj}^{+}$, we know $L_{mj}\setminus W_{mj}^{+}$
is closed in $L_{mj}$.
We next show that $L_{mj}\setminus W_{mj}^{+}$ is of Hausdorff codimension at least $2$.
Consider the restriction of the map $E$ on $L_{mj}\setminus W_{mj}^{+}$ given by
$$
f:=(p,t)\in L_{mj}\setminus W_{mj}^{+}\mapsto E(p,t)\in {\mathcal F}.
$$
Note that $f$ is a finite-to-one map. This is because for any $z\in {\mathcal F}$, the height
function $\ell_{z}$ is a perfect Morse-Bott function, and each $p$, for which some point of the form
$(p,t)$ is in $f^{-1}(z)$, is a nondegenerate critical point of index $m$ of
$\ell_z$ on the taut
hypersurface $M$.  Such a nondegenerate critical point
is a critical submanifold. Hence the total number of such points is no more than $b_{m}$,
the $m$-th Betti number of $M$.

The image of $f$ is a closed subset of ${\mathcal F}$ of Hausdorff codimension
at least $2$ by Lemma~\ref{Hausdorff}. $E$ is
a local diffeomorphism
at $(p,t)\in L_{mj}$ when it is restricted to $L_{mj}$. Hence, there is a
neighborhood $Y_{(p,t)}\subset L_{mj}$ around each
$(p,t)\in L_{mj}\setminus W_{mj}^{+}$
such that $f$ restricted to $Y_{(p,t)}\cap (L_{mj}\setminus W_{mj}^{+})$ is a
homeomorphism into
${\mathcal F}$. It follows that the Hausdorff $(n-2)$-measure of
$L_{mj}\setminus W_{mj}^{+}$
is at most $b_{m}$ times that of ${\mathcal F}$, which is null in the ambient sphere
$S^{n}$. So,
$L_{mj}\setminus W_{mj}^{+}$ is of Hausdorff codimension at least 2 in $L_{mj}$.
In particular, $L_{mj}\setminus W_{mj}^{+}$
does not disconnect $L_{mj}$.

We now apply Remark~\ref{fractal} with
$S=L_{mj}\setminus W_{mj}^{+}$, $M=L_{mj}$ and $X = M \setminus S = W_{mj}^{+}$.  
Since $X$ has a finite number of ends, it follows that $S$ has finitely
many components all of which are open in $S$.
Now, suppose $M$ has infinitely many ends. Then
there is a decreasing sequence of open sets $X_{1},X_{2},\ldots$, in $X$
such that $X\setminus X_{j}$ is compact, $\cap_{j} \overline{X}_{j}=\emptyset$ and the number of unbounded
connected components of $X_{j}$
increasingly diverges to infinity. As $S$ has finitely many components $S_{1},\ldots,S_{t}$,
all of which are of Hausdorff codimension at least $2$, $S_{1},\ldots,S_{t}$
cannot disconnect the connected components of $X_{j}$. We see $X_{j}\setminus S$ form
a decreasing sequence whose disconnected components diverge to infinity, so that
$W_{mj}^{+}=M\setminus S$ will have infinitely many ends. This contradiction establishes
that $L_{mj}$ has a finite number of ends. Since $L_{mj}$ is a trivial
line bundle over $U_{mj}^{+}$, so that $L_{mj}\simeq U_{mj}^{+}\times \R$, it follows that
$U_{mj}^{+}$ has a finite number of ends.
\end{proof}

\begin{corollary}\label{finite1} $W_{mj}^{+}$ is disjoint from $W_{ml}^{+}$
if and only if $U_{mj}^{+}$
is disjoint from $U_{ml}^{+}$. In particular, $U_{m}^{+}$ has finitely many ends.
\end{corollary}

\begin{proof} The backward direction is clear. To prove the forward direction,
suppose $p\in U_{mj}^{+}\cap U_{ml}^{+}$. Let $X:=U_{mj}^{+}\cup U_{ml}^{+}$ and let
$Y:=W_{mj}^{+}\cup W_{ml}^{+}$. Since $X$ is open and connected, we can form the
(trivial) line bundle $L$ over $X$, where $L$ consists of $(p,t),p\in X,$ such
that the distance function $d_{q}$, for $q=E(p,t)$, has a nondegenerate critical point of index $m$ at $p$. 
Then $L=L_{mj}\cup L_{ml}$ and
$Y\subset L$ by construction. However, the analysis in
Lemma~\ref{finite-ends}
shows that $L \setminus Y$ is of Hausdorff codimension at least $2$; therefore, $Y$ is
connected, and so $W_{mj}^{+}$
and $W_{ml}^{+}$ cannot be disjoint.

Now that $U_{m}^{+}$ is the finite disjoint union of all $U_{mj}^{+}$, each of
which has
finitely many ends, it follows that the same is true for $U_{m}^{+}$.
\end{proof}

\begin{corollary}\label{*set1} $U_{m}^{+}$ is dense in
$(U_{m}^{*})^{+}$. In particular, $(U_{m}^{*})^{+}$ has a finite number of connected
components.
\end{corollary}

\begin{proof} Suppose that $U_{m}^{+}$ is not dense in
$(U_{m}^{*})^{+}$. Then there is an open set $X\subset (U_{m}^{*})^{+}$ that is disjoint from $U_{m}^{+}$.
Similar to Corollary~\ref{finite1}, let us introduce the line bundle $L$
over $X$, which consists of all
points $(p,t),t>0,$ where $p\in X$ and $d_{q}$, for $q=E(p,t)$, 
has a nondegenerate critical point of index $m$ at $p$.

Every point $(p,t)\in L$ is mapped to ${\mathcal F}$ via the map $E$, since $p \in (U_{m}^{*})^{+} \setminus U_{m}^{+}$.
As we see in Lemma~\ref{finite-ends}, the map $f$ of Lemma~\ref{finite-ends} on $L$ is finite-to-one and regular
at every point of $L$.  Thus, the image of $f$ contains an open set of dimension $n = \dim L$.  This contradicts the fact
that the Hausdorff codimension of ${\mathcal F}$ is at least two. Thus, $U_{m}^{+}$ is dense in
$(U_{m}^{*})^{+}$.

In particular, $(U_{m}^{*})^{+}$ also has a finite number of connected components,
because any of its
components will contain at least one component of $U_{m}^{+}$.
\end{proof}

Recall that we say a point $p$ is a good point in the taut hypersurface $M$
if the multiplicities of the principal curvatures are locally constant around $p$.
We denote the set of good points by ${\mathcal G}$. We know ${\mathcal G}$
is open and dense in $M$. 

Our convention is that we label the principal curvature functions by
$\lambda_{1}\leq\lambda_{2}\leq \cdots \leq\lambda_{n-1}$. When we say
$(m_{1},m_{2},\ldots,m_{g})$
is a given sequence of principal multiplicities, we mean that
\begin{displaymath}
\lambda_1 = \cdots = \lambda_{m_1}, \quad \lambda_{{m_1 + 1}} = \cdots = \lambda_{m_1 + m_2}, \ldots,
\end{displaymath}
so that $m_i$ is the multiplicity of the $i^{\rm th}$ largest distinct principal curvature.

\begin{lemma}\label{approx} There is a dense open subset of ${\mathcal G}$
such that for any $p$ in the subset,
there is
a point $q\in {\mathcal Z}$
between any two focal points on the normal exponential circle $E(p,t),t>0.$
\end{lemma}

\begin{proof}
Let $X$ be a connected component of ${\mathcal G}$ assuming
multiplicities $(m_{1},m_{2},\ldots,m_{\mu})$.
Let
$\sigma_{s}=m_{1}+\cdots+m_{s},$ for $1 \leq s \leq \mu$. Then $X\subset (U_{\sigma_{s}}^{*})^{+}$.
Then Corollary~\ref{*set1} implies that $X\cap U_{\sigma_{s}}^{+}$
is dense and open in $X$.
\end{proof}

\begin{definition} We let ${\mathcal G}^{\circ}$ be the largest dense open
subset of ${\mathcal G}$ over which Lemma~\ref{approx} is true.
\end{definition}

\begin{lemma}\label{pinpoint} 
Let $g$ be the maximum number of distinct principal curvatures on $M$.
Let ${\mathcal M}:=(m_{1},\ldots, m_g)$ be a given maximal sequence of
multiplicities. Let $O_{\mathcal M}\subset {\mathcal G}^{\circ}$ be the (open)
subset of ${\mathcal G}^{\circ}$
attaining these multiplicities. Let 
$\sigma_{j}=m_{1}+\cdots+m_{j}$, for $1 \leq j \leq g$. Then $O_{\mathcal M}$ is
the intersection of all
$U_{\sigma_{j}}^{+}$.
\end{lemma}

\begin{proof} The intersection consists of all points $p$ such that
there is some $0<t_{j}<\pi$ for which $p$ is of index $\sigma_{j}$ with
respect to the distance function $d_{q}$, with $q=E(p,t_{j}),$ for each $j$.
So,
$p$ must have $g$ principal curvatures with multiplicities
$m_{1},m_{2},\ldots,m_g$.
The maximality of $g$ implies that $m_{j}$ cannot be broken further into smaller
multiplicities. Moreover, since $p\in U_{\sigma_{j}}^{+}$, there is a $t$ for
which $E(p,t)$ between the two appropriate focal points is in ${\mathcal Z}$. We
obtain $p\in {\mathcal G}^{\circ}$. So, $\cap U_{\sigma_{j}}^{+}\subset {\mathcal O}_{M}$.

Conversely, since ${\mathcal O}_{M}\subset {\mathcal G}^{\circ}$, it follows
by the definition of $U_{\sigma_{j}}^{+}$ that
${\mathcal O}_{M}\subset\cap U_{\sigma_{j}}^{+}$.
\end{proof}

We now proceed to handle the case of a taut submanifold of dimension 3 or 4.  First we need the following lemma 
in homology that will be needed in our proof.

\begin{lemma}
\label{homology} 
Let $B_1 \subset B_2$ be two closed disks. Then the $k^{\rm th}$ new
topology of the taut submanifold $M$ added between $B_1 \cap M$ and $B_2 \cap M$ is the relative homology
\begin{displaymath}
H_k(B_2 \cap M, B_1 \cap M).
\end{displaymath}
\end{lemma}

\begin{proof}
This follows from the exact sequence
\begin{eqnarray}
\rightarrow H_k(B_1\cap M) & \stackrel{i}{\rightarrow} & H_k(B_2\cap M) \rightarrow  H_k(B_2\cap M, B_1\cap M)\nonumber \\
 & \stackrel{j}{\rightarrow} &  H_{k-1}(B_1\cap M) \stackrel{i}{\rightarrow} H_{k-1}(B_2\cap M) \rightarrow \nonumber
\end{eqnarray}
and the fact that $i$ is injective by the injectivity property of tautness,
so that $j$ is surjective. It follows that
\begin{displaymath}
H_k(B_2\cap M) = H_k(B_1\cap M) \oplus H_k(B_2\cap M, B_1\cap M),
\end{displaymath}
and hence the conclusion.
\end{proof}

\begin{corollary}
\label{cor-homology}
If $B_1\subset B_2\subset B_3$ are three closed disks, then
\begin{displaymath}
H_k(B_3\cap M, B_1\cap M) = H_k(B_3\cap M, B_2\cap M) \oplus H_k(B_2\cap M, B_1\cap M).
\end{displaymath}
Thus, $H_k(B_3\cap M, B_1\cap M)$ is surjective to both
$H_k(B_3\cap M, B_2\cap M)$ and $H_k(B_2\cap M, B_1\cap M)$.
\end{corollary}
\begin{proof}
Given three abelian groups $A \rightarrow B \rightarrow C$, where each arrow is an
embedding, then $C/B = (C/A)/(B/A)$.
\end{proof}

\begin{lemma}
\label{cobordism} 
Suppose $\dim M=4$.\\
{\rm (a)}  If the maximum number of principal curvatures
is $\geq 3$, then multiplicities $(1,3),(3,1)$ and $(2,2)$ cannot exist on open sets. \\
{\rm (b)}If the maximum number of principal curvatures is $4$, then the points with
multiplicities $(1,1,2),(1,2,1)$ and $(2,1,1)$, at which the number of principal
curvatures is not locally constant, cannot be approached by a sequence coming from an open set of
points of the same multiplicities.
\end{lemma}

\begin{proof}
(a) Suppose the multiplicities are $(1,3)$ on an open set $O$. Let $x \in O$ and let $p$ be a boundary point of $O$.
Let $c(t)$ be a smooth curve with $c(0) = x$ and $c(1) = p$.  We can assume that $c(t) \in O$ for $0 \leq t < 1$.
Otherwise, let $t_0$ be the first value of $t$ such that $q = c(t_0)$ 
is not in $O$, and replace $p$ by $q$.
 
Since $p$ is a boundary point of $O$, there must be a sequence of points with
multiplicities $(1,1,2)$, $(1,2,1)$ or $(1,1,1,1)$ that converges to $p$.  
The multiplicities must remain $(1,3)$ at $p$; otherwise, the list of multiplicities would drop to the single multiplicity
$(4)$, and this is impossible, because a taut hypersurface with a single umbilic point
must be a totally umbilic sphere.  We will handle the case where a sequence of points with 
multiplicities $(1,1,2)$ approaches $p$.  The other cases are very similar.

At each point $c(t), 0 \leq t < 1$, the curvature surface $S(t)$ corresponding
to the principal curvature $\lambda$ of multiplicity 3 is a 3-dimensional metric sphere which is the intersection of the
corresponding 4-dimensional curvature sphere $\Sigma (t)$ with a 4-dimensional plane $P(t)$. As we take the limit as $t$ 
approaches $1$, these 4-dimensional planes $P(t)$ approach a limiting 4-plane $P(1)$ that intersects the
4-dimensional curvature sphere $\Sigma (1)$ in a 3-dimensional sphere $S(1)$, which must be the 3-dimensional
curvature surface corresponding to $\lambda$ at $p$.  On the other hand, consider a sequence of points $\{y_i\}$
where the multiplicities are $(1,1,2)$ that approaches $p$.  For each $y_i$, the 2-dimensional
curvature surface $C_i$ through $y_i$ corresponding to the principal curvature of multiplicity 2 is a topset for
the spherical distance function centered at the corresponding focal point.  Thus, by tautness, $C_i$ represents a
nontrivial 2-dimensional homology class in $M$.  As $y_i$ approaches $p$, these $C_i$ approach a 2-cycle $C$ in
the curvature surface $S(1)$ that is nontrivial in the 2-dimensional homology of $S(1)$.  This is a contradiction,
since $S(1)$ is a 3-sphere and has trivial 2-dimensional homology.
A similar proof shows that multiplicities $(3,1)$ cannot exist on an open set.

Next suppose that the multiplicities are $(2,2)$ on an open set $O$. As in the argument above,
let $x \in O$ and let $p$ be a boundary point of $O$.
Let $c(t)$ be a smooth curve with $c(0) = x$, $c(1) = p$ and $c(t) \in O$ for $0 \leq t < 1$.  
Then the multiplicities must be $(2,2)$ at $p$, because neither of the two multiplicities can drop
to 1 at the limit point $p$ and the single multiplicity $(4)$ is impossible.  
Since $p$ is a boundary point of $O$, there must be a sequence of points with
multiplicities $(1,1,2)$, $(2,1,1)$ or $(1,1,1,1)$ that converges to $p$. We will handle the case where the
multiplicities are $(1,1,2)$, and the others are handled similarly.

Then as in the argument above, 
at each point $c(t), 0 \leq t < 1$, the curvature surface $S(t)$ corresponding
to the first principal curvature $\lambda$ of multiplicity 2 is a 2-dimensional metric sphere which is the 
intersection of the
corresponding 4-dimensional curvature sphere $\Sigma (t)$ with a 3-dimensional plane $P(t)$. As we take the limit as $t$ 
approaches $1$, these 3-dimensional planes $P(t)$ approach a limiting 3-plane $P(1)$ that intersects the
4-dimensional curvature sphere $\Sigma (1)$ in a 2-dimensional sphere $S(1)$, which must be the 2-dimensional
curvature surface corresponding to $\lambda$ at $p$. On the other hand, there is a sequence $\{y_i\}$
approaching $p$ such that the multiplicities are $(1,1,2)$ at $y_i$.  For each $y_i$, the 1-dimensional
curvature surface $C_i$ through $y_i$ corresponding to the first principal curvature of multiplicity 1 is a topset for
the spherical distance function centered at the corresponding focal point. Thus, $C_i$ represents a
nontrivial 1-dimensional homology class in $M$.  As $y_i$ approaches $p$, these $C_i$ approach a 1-cycle $C$ in
the curvature surface $S(1)$ that is nontrivial in the 1-dimensional homology of $S(1)$.  This is a contradiction,
since $S(1)$ is a 2-sphere and has trivial 1-dimensional homology.
\\
(b) In the case of multiplicities $(1,1,2)$ or $(2,1,1)$, the same type of argument takes care of this statement, since 
one can produce a point $p$ with the given multiplicities
that is also a limit of a sequence of points with multiplicities $(1,1,1,1)$.
Then the nondegenerate critical manifold at $p$ corresponding to the principal curvature of multiplicity two
is a metric 2-sphere $S$.  However, since $p$ can be approached by a sequence of points with multiplicities
$(1,1,1,1)$, one can produce a top 1-cycle $C$ in $S$ that is nontrivial in the 1-dimensional homology of $S$, which
is a contradiction. 

In the case of multiplicities $(1,2,1)$, we need to modify the argument slightly. As above,
the nondegenerate critical manifold at $p$ corresponding to the principal curvature of multiplicity two
is a metric 2-sphere $S$, and we can also approach $p$ by a sequence of points with multiplicities $(1,1,1,1)$.
However, in this case, the 1-dimensional circles approaching $p$ corresponding to the second and third multiplicities
are not top sets.  We can see that they still contribute in a nontrivial way to homology by the following argument. 

At the point $p$, the multiplicities are $(1,2,1)$, and $p$ is a limit point of the open connected set $O$ on
which the multiplicities are $(1,2,1)$.  Let $q$ be the second focal point of $M$ at $p$ corresponding to the
principal curvature of multiplicity two.  As above, we can also approach $p$ by a sequence of points in the open set $U$ with
multiplicities $(1,1,1,1)$.  Thus near $p$, we can find points $x \in O$ and $y \in U$, with corresponding second focal
points $u$ and $v$ near $q$ such that 
the height functions $\ell_u (x) = a$, $\ell_v (y) = b$, and positive numbers $\epsilon$ and
$\delta$ such that
\begin{equation}
\label{Morse}
M_{a-\epsilon}(\ell_u) \subset M_{b -\delta}^{-} (\ell_v) \subset M_b(\ell_v), \quad 
M_a (\ell_u) \subset M_{b + \delta}^{-}(\ell_v) \subset M_{a+\epsilon}(\ell_u).
\end{equation}
Furthermore, using a genericity argument, we can assume that the points $x$ and $y$ are chosen so that the critical
submanifolds of the corresponding height functions $\ell_u$ and $\ell_v$ are at distinct levels, and we may choose
$\epsilon$ and $\delta$ sufficiently small so that $a$ (respectively $b$) is the only critical value between
$a - \epsilon$ and $a + \epsilon$ (respectively, between $b - \delta$ and $b + \delta$). 

Lemma~\ref{homology} says that the new $k^{\rm th}$ homology between the levels $a-\epsilon$ and $a+\epsilon$ of $\ell_u$ is
\begin{equation}
\label{eq:h1}
H_k (M_{a+\epsilon} (\ell_u), M_{a-\epsilon} (\ell_u)).
\end{equation}
Corollary~\ref{cor-homology} says that the group in equation \ref{eq:h1} is surjective to
\begin{equation}
\label{eq:h2}
H_k (M_{a+\epsilon} (\ell_u), M_{b - \delta} (\ell_v)),
\end{equation}
where in Lemma~\ref{homology}, we take
\begin{displaymath}
B_1 \cap M = M_{a - \epsilon} (\ell_u),  B_2 \cap M = M_{b - \delta} (\ell_v),  B_3 \cap M = M_{a + \epsilon} (\ell_u).
\end{displaymath}
The same lemma says that the group in equation \ref{eq:h2} is surjective to
\begin{equation}
\label{eq:h3}
H_k (M_{b + \delta} (\ell_v), M_{b - \delta} (\ell_v)),
\end{equation}
where
\begin{displaymath}
B_1 \cap M = M_{b - \delta} (\ell_v),  B_2 \cap M = M_{b + \delta} (\ell_v),  B_3 \cap M = M_{a + \epsilon} (\ell_u).
\end{displaymath}
It follows from these considerations that the group in equation \ref{eq:h1} is surjective to the group in
equation \ref{eq:h3}.  However, by Morse-Bott critical point theory (see \cite[Theorem 20.2, p. 239]{DuFoNoIII}), we have
\begin{equation}
\label{eq:h4}
H_k (M_{a + \epsilon} (\ell_u), M_{a - \epsilon} (\ell_u)) = H_{k-\mu} (W),
\end{equation}
where $W$ is the critical manifold of $\ell_u$ at $x$ and $\mu$ is the index at $x$.  Similarly, we have
\begin{equation}
\label{eq:h5}
H_k (M_{b + \delta} (\ell_v), M_{b - \delta} (\ell_v)) = H_{k-\nu} (V),
\end{equation}
where $V$ is the critical submanifold of $\ell_v$ at $y$ and $\nu$ is the index at $y$.

Now the critical submanifold $W$ at $x$ is a 2-sphere of index $\mu = 1$, whereas the critical submanifold
$V$ at $y$ is a circle of index $\nu = 1$.  Using $k = 2$ in equations \ref{eq:h4} and \ref{eq:h5}, we see that 
$0 = H_1(W,\Z_2)$ is surjective to $H_1 (V,\Z_2) = \Z_2.$  This is a contradiction.
\end{proof}

\begin{theorem}
\label{dim3-4}
Let $M$ be a compact taut hypersurface in $S^{n}$ that is not a hypersphere.\\
{\rm (I)} If $\dim M=3$, then $M$ is algebraic.\\
{\rm (II)} If $\dim M=4$, then $M$ is algebraic.\\
\end{theorem}

\begin{proof}
(I) We first establish that generically, i.e., on a dense open set,
there is only one sequence of multiplicities.
For $\dim M=3$, if no points assume multiplicities $(1,1,1)$, then we have 
the proper Dupin case of two distinct principal curvatures at every point.  This follows from the fact that 
no umbilic points can exist, since $M$ is not a hypersphere. In this case, $M$ is algebraic by Corollary~\ref{Main}.

In the case where there is an open set on
which the multiplicities are $(1,1,1)$,
the same type of argument used in the multiplicities $(1,3)$ case in Lemma~\ref{cobordism}
shows that the multiplicities $(2,1)$ and $(1,2)$ cannot exist on open sets.  Thus, the multiplicities
are $(1,1,1)$ on the dense open subset ${\mathcal G}$ of $M$.

Now, consider $A:=(U_{1}^{*})^{+}$ and $B:=(U_{1}^{*})^{-}$. The set $A$ consists of points
of multiplicities $(1,1,1)$ or $(1,2)$ and $B$ with multiplicities $(1,1,1)$
or $(2,1)$.  The sets $A$ and $B$ both have finitely many connected components, and their union is all of $M$.
Let $W:=A\cap B$ be the set on which the multiplicities are $(1,1,1)$. The Mayer-Vietoris sequence 
$$
\leftarrow H^{1}(M)\leftarrow H^{0}(W)\leftarrow H^{0}(A)\oplus H^{0}(B)
\leftarrow H^{0}(M)\leftarrow 0
$$
gives that $W = {\mathcal G}$ has finitely many
connected components. Hence, the local finiteness property holds on $M$, and $M$ is algebraic by Theorem~\ref{Kuiper}.

(II)(a) Let $\dim M = 4$.  First consider the case where the maximum number of distinct principal curvatures is two.
Since $M$ is not a hypersphere, there cannot be an umbilic point, and therefore the number of distinct principal
curvatures must be two at all points.  Thus, $M$ is proper Dupin, and $M$ is algebraic by Corollary~\ref{Main}.

Next consider the case when the maximum number of distinct principal curvatures is three.
Note, in general, that
the set assuming multiplicities $(1,1,2)$, $(1,2,1)$ or $(2,1,1)$
is automatically open by the maximality of the
number of distinct principal curvatures. As a consequence, only the multiplicities $(1,1,2),(1,2,1),(2,1,1)$ can possibly
exist on open sets by Lemma~\ref{cobordism}.

Consider $A:=(U_{1}^{*})^{+},B:=(U_{2}^{*})^{+}$ and $C:=(U_{3}^{*})^{+}$. The set
$A$ consists of points with multiplicities $(1,1,2),(1,2,1),(1,3)$, $B$ with
multiplicities $(1,1,2),(2,1,1),(2,2)$, and $C$ with multiplicities
$(1,2,1)$, $(2,1,1)$, $(3,1)$. Each of $A,B,C$ has finitely many
connected components, and $M$ is the union of these three sets.

Consider $A$ and $D:=B\cup C$. Being a union, $D$ also has a finite number of connected components.
Now $M=A\cup D$ and $W:=A\cap D=U\cup V$, where $U$ is the open set of
points with multiplicities $(1,1,2)$ and $V$ is the open set of points
with multiplicities $(1,2,1)$. Hence, the Mayer-Vietoris sequence applied to
$A$ and $D$
implies that $W$ has finitely many connected components. As the open sets
$U$ and $V$ are disjoint, it follows that each of $U$ and $V$ has finitely
many connected 
components. A similar consideration establishes that the set of points
with multiplicities
$(2,1,1)$ also has a finite number of components by considering $B$ and $A\cup C$. Hence, the set
${\mathcal G}$ has finitely many
connected components.
Therefore, the local finiteness property holds on $M$, and $M$ is algebraic by Theorem~\ref{Kuiper}.

(II)(b) Assume now that the maximum number of multiplicities is $4$. By
Lemma~\ref{cobordism},
multiplicities $(1,3),(3,1)$ and $(2,2)$ cannot exist on open sets. Meanwhile,
the set of points with multiplicities $(1,1,1,1)$ is open.

Since by Lemma~\ref{cobordism}, a point with multiplicities $(1,1,2)$
(or $(1,2,1)$ or $(2,1,1)$) at which the number of principal curvatures is not
locally constant cannot be approached by a sequence 
coming from an open set of points of the same multiplicities, we see that
such points must be entirely surrounded by points with multiplicities
$(1,1,1,1)$.
As a result, an open set of points with multiplicities $(1,1,2)$
will approach boundary points with multiplicities either $(1,3)$ or $(2,2)$.

The set ${\mathcal G}$ of good points are those with multiplicities $(1,1,1,1)$
and the points with multiplicities $(1,1,2),(1,2,1)$, or $(2,1,1)$ that exist on open
sets. ${\mathcal G}$ is open and dense in $M$. We let $S$ be the subset of
${\mathcal G}^{c}$ consisting of points
of multiplicities $(1,3),(3,1)$ or $(2,2)$. We wish to establish that the local
finiteness property is true.

Firstly, $S$ is closed in $M$. This is because any converging sequence of
points of the indicated multiplicities must maintain the same type of
multiplicities.

We next show that $M\setminus S$ has finitely many connected components. It comes down to
showing that $A:=(U_{1}^{*})^{+}$ (respectively, $B:=(U_{2}^{*})^{+}$ or
$C:=(U_{3}^{*})^{+}$),
with points of multiplicities $(1,3)$ (respectively, $(2,2)$ or $(3,1)$)
removed, has only finitely many connected components; for then the union of
the three resulting sets is exactly $M \setminus S$, which, being a union, must have
finitely many connected components. However, as before, this follows from the
Mayer-Vietoris sequence applied to the open covers $A$ and
$B\cup C$ of $M$, etc. Note that $A$ is composed of points with multiplicities
$(1,1,1,1),(1,1,2),(1,2,1)$ or $(1,3)$, $B$ with multiplicities
$(1,1,1,1),(1,1,2),(2,1,1)$ or $(2,2)$, and $C$ with multiplicities
$(1,1,1,1),(1,2,1),(2,1,1)$
or $(3,1)$. Hence, $A\cap(B\cup C)$ is exactly $A$ with points of
multiplicities $(1,3)$ removed.

Lastly, we verify that each point in ${\mathcal G}^{c} \setminus S$ has a
open neighborhood
$W$ in $M$ such that $W\cap{\mathcal G}$ contains finitely many connected
open sets whose union is dense in $W$.
















As mentioned earlier, ${\mathcal G}^{c}\setminus S$ is completely surrounded
by points of multiplicities $(1,1,1,1)$. For each point
$p\in {\mathcal G}^{c}\setminus S$ of multiplicities $(1,1,2)$, there is a small
neighborhood of it which contains no points of multiplicities $(1,2,1)$ or
$(2,1,1)$, etc.
In other words, The sets $T_1,T_2,T_3$ of points of multiplicities $(1,1,2),
(1,2,1),(2,1,1),$
respectively, are contained in disjoint open sets $O_1,O_2,O_3$, respectively,
where $O_i\setminus T_i,1\leq i\leq 3,$ consists of only points of multiplicities
$(1,1,1,1)$.

Let $X\subset (U_3)^{+}$ be the subset of points of multiplicities $(1,1,1,1)$.
$X$ is dense in the set of points of multiplicities $(1,1,1,1)$ by
Lemmas~\ref{approx} and~\ref{pinpoint}. It follows that
$T_1$ is contained in the end sets
of the finitely many end components $E_1,\ldots,E_s$ of
$U_3^{+}$ on which the multiplicities remain $(1,1,1,1)$. Let $W$ be the
union of $O_1$ above and
$E_1,\ldots,E_s$. ($O_1$ is the union of neighborhoods
of points of multiplicities $(1,1,2)$. We make sure each of these
neighborhoods is
so small that its intersection with $U_3^{+}$ is contained in the end
components
$E_1,\ldots,E_s$.) Then $W\cap {\mathcal G}$ consists of only points of
multiplicities
$(1,1,1,1)$, and moreover contains open sets $E_1,\ldots,E_s$ whose union is
dense in $W$.









Similarly, the same conclusion is true for points of multiplicities
$(2,1,1)$ and $(1,2,1)$ in ${\mathcal G}^{c}$ with $U_3^{+}$ replaced by
$U_1^{+}$ and $U_2^{+}$, respectively. Hence, $M$ is algebraic.
\end{proof}

\vspace{3mm}
\begin{remark} In fact, the proof of (I) and (II)(a) are the same as that of
(II)(b) in a hidden way. For (II)(a), ${\mathcal G}^{c}$ is the set of points
with multiplicities $(1,3),(3,1)$ or $(2,2)$, which coincides with $S$, so that the set
${\mathcal G}^{c} \setminus S=\emptyset$, and the local finiteness condition is
automatically satisfied on ${\mathcal G}^{c} \setminus S$.
Similarly for (I), we have ${\mathcal G}^{c}=S$, which is the set of points with
multiplicities $(1,2)$ or $(2,1)$.
\end{remark}

Recall that an embedding $f:M \rightarrow \R^n$ is said to be {\em substantial} if the image of $f$
does not lie in any affine hyperplane in $\R^n$.

\begin{theorem}
\label{dim<5}
Let $M$ be a compact taut $m$-dimensional submanifold in $S^n \subset \R^{n+1}$.  If $m \leq 4$, then $M$ is algebraic.
\end{theorem}

\begin{proof}
Banchoff \cite{Ban} showed that a taut compact 1-dimensional submanifold of $S^n$ must be a metric circle in $S^n$,
which is certainly algebraic. In the same paper, 
he also showed that if $M$ is a taut compact 2-dimensional surface
substantially embedded in $S^n$, then $M$ is a metric 2-sphere,
a cyclide of Dupin in $S^3$, or a spherical Veronese surface $V \subset S^4$.  All of these 
surfaces are algebraic.

Next let $M$ be a compact taut 3-dimensional submanifold of codimension $k+1$ in $S^n \subset \R^{n+1}$. Without loss
of generality, we may assume that $M$ is substantially embedded in $\R^{n+1}$.  We want to consider the tube
$M_\epsilon$ of radius $\epsilon > 0$ over $M$ in $S^n$, where $\epsilon$ is sufficiently small so that $M_\epsilon$
is an embedded hypersurface in $S^n$.  By a theorem of Pinkall \cite{Pinkall2}, we know that $M_\epsilon$
is a taut hypersurface in $S^n$.  We can parametrize $M_\epsilon$ by the map 
$\phi_\epsilon: UN(M) \rightarrow S^n$, where $UN(M)$ is the unit normal bundle of $M$ in $S^n$, given by
\begin{displaymath}
\phi_\epsilon (x,\xi) = \cos \epsilon x + \sin \epsilon \xi,
\end{displaymath}
where $\xi$ is a unit normal vector to $M$ at $x$.  Then
\begin{displaymath}
\eta = - \sin \epsilon x + \cos \epsilon \xi,
\end{displaymath}
is a unit normal vector to the tube $M_\epsilon$ at the point $y = \phi_\epsilon (x,\xi)$.  In 
\cite[pp. 131-132]{CecilRyan} a formula is given for the principal curvatures of the shape operator
$A_\eta$ of $M_\epsilon$ at $y$ in terms of the shape operator $A_\xi$ of $M$ at $x$.  There it is shown that
$A_\eta$ has a principal curvature $\mu_0 = \cot (-\epsilon)$ of multiplicity $k$ and principal curvatures
\begin{displaymath}
\mu_i = \cot (\theta_i - \epsilon), \quad 1 \leq i \leq 3,
\end{displaymath}
where $\lambda _i = \cot \theta_i, 1 \leq i \leq 3$, for
$0 < \theta_3 < \theta_2 < \theta_1 < \pi$, are
the principal curvatures of $A_\xi$ at $x$.  Thus each $\mu_i, 1 \leq i \leq 3$, has the same multiplicity as
the corresponding $\lambda_i$.  Therefore, on $M_\epsilon$ the multiplicities of the principal curvatures
must take the form $(k,1,1,1), (k,1,2)$ or $(k,2,1)$.  Note that $(k,3)$ is not possible, for then $A_\xi$
would have a principal curvature of multiplicity 3 at $x$, i.e., all the principal curvatures of $A_\xi$ would
be equal at $x$.  Then, if $p$ is the first focal point along the normal geodesic in $S^n$ from $x$
in the direction $\xi$, tautness implies that the height function $\ell_p$ must have both an absolute
maximum and an absolute minimum at $x$ (see \cite[Lemma 1.24, p. 122]{CecilRyan}).  Thus all of $M$ lies in
in hyperplane in $\R^{n+1}$ orthogonal to $p$.  This contradicts the assumption that $M$ is substantial
in $\R^{n+1}$, and so the multiplicities $(k,3)$ cannot occur on $M_\epsilon$.

As in the 3-dimensional hypersurface case in Theorem~\ref{dim3-4}, multiplicities $(k,1,2)$, respectively, $(k,2,1)$, cannot
exist on open sets in $M_\epsilon$, unless the multiplicities have the constant values $(1,2)$, respectively, $(2,1)$,
on the unit normal bundle $UN(M)$ of $M$.  In that case, the tube $M_\epsilon$ is a proper Dupin hypersurface,
and so $M_\epsilon$ and $M$ are algebraic by Theorem~\ref{Kuiper}.
Hence, the only remaining case is when multiplicities have the values $(1,1,1)$ on a dense open subset of $UN(M)$.

Then another application of the Mayer-Vietoris sequence applied to the sets
$(U_{1}^{*})^{+}$ and $(U_{2}^{*})^{+}$ for the tube $M_\epsilon$ shows that the generic multiplicities $(k,1,1,1)$
exist on only finitely many connected components.  Thus, $M_\epsilon$ satisfies the local finiteness property 
and is algebraic, and so $M$ is algebraic by Theorem~\ref{Kuiper}.

Suppose now that $M$ is a compact taut 4-dimensional submanifold substantially embedded in $S^n \subset \R^{n+1}$
with codimension $k+1$ in $S^n$.  By the same construction as above, the possible multiplicities on $M_\epsilon$ are
\begin{displaymath}
(k,1,1,1,1), (k,1,1,2), (k,1,2,1), (k,2,1,1), (k,1,3), (k,3,1), (k,2,2).
\end{displaymath}  
Then the same arguments used to prove
Lemma~\ref{cobordism} and Theorem~\ref{dim3-4} (II) can be applied here by just adjoining the multiplicity $k$ at the 
beginning, and thus we conclude that $M_\epsilon$ is algebraic, and so $M$ is algebraic.

\end{proof}

\bibliography{Bibliography}
\bibliographystyle{plain}

\end{document}